\documentclass[11pt,a4paper]{article}

\usepackage[round,authoryear]{natbib}
\usepackage[hidelinks]{hyperref}

\usepackage[utf8]{inputenc}
\usepackage[T1]{fontenc}
\usepackage{lmodern}
\usepackage{amsmath, amssymb, amsthm}
\usepackage{graphicx}
\usepackage{caption}
\usepackage{subcaption}
\usepackage{geometry}
\usepackage{float}
\usepackage{tikz}
\usetikzlibrary{calc}
\usetikzlibrary{fit}
\usetikzlibrary{matrix, positioning, arrows.meta, shapes.geometric}
\usepackage{xcolor}
\usepackage{indentfirst}
\usepackage{authblk}

\theoremstyle{definition}
\newtheorem{definition}{Definition}[section]

\newtheorem{remark}{Remark}

\title{A system-centered probabilistic formalism linking multiple equilibria and biodiversity in ecological and evolutionary models}

\author[1,*]{Hiba Nassor}
\author[1]{Hermine Biermé}
\author[2]{Elisabeth Herniou}
\author[1]{Sten Madec}

\affil[1]{Institut Denis Poisson, CNRS, IDP UMR 7013, Université de Tours, Université d'Orléans, Tours, France}
\affil[2]{Institut de Recherche sur la Biologie de l'Insecte, UMR CNRS 7261, Université de Tours, Tours, France}
\affil[*]{Corresponding author}

\date{}

\begin{document}

\maketitle

\begin{center}
\small
Corresponding author: hiba.nassor@univ-tours.fr
\end{center}

\vspace{0.5cm}

\begin{abstract}

Ecological and evolutionary communities do not always settle into a single predictable configuration. Models such as Lotka-Volterra or replicator dynamics may predict several admissible states. The multiplicity of such states is central to understanding robustness, alternative community configurations, coexistence, and shifts in biodiversity. Yet, standard practices are often state-centered: they either focus on the properties of typical equilibria or pool states across systems, losing information about how many states can arise in a single system, including cases where no state exists.

Here, we introduce a system-centered probabilistic formalism that captures this hidden structure and could be applied across a wide range of ecological and evolutionary dynamical models, drawing inspiration from hurdle and zero-inflated models. Each ecological system is viewed as generating a distribution over its possible outcomes. This allows us to jointly quantify the occurrence of a given state type, its multiplicity within each system, and the number of species co-existing within these states.

The formalism is generic and can be applied to different families of dynamical models. To illustrate its scope, we consider two case studies: the generalized Lotka-Volterra and replicator dynamics with random and structured parameters. This unified perspective reveals patterns that remain invisible when equilibria are analyzed individually or pooled across systems.

By turning the state space of a model into interpretable probabilistic quantities, our formalism offers a new way for studying the robustness of community structure, the likelihood of alternative outcomes under small changes in initial conditions, and the potential for shifts between states of low and high diversity.

\end{abstract}


\vspace{0.5cm}
\noindent\rule{\textwidth}{0.4pt}
\vspace{0.3cm}

\noindent\textbf{\small Funding} \\
\small This project has received financial support from the CNRS through the MITI interdisciplinary programs.

\vspace{0.4cm}

\noindent\textbf{\small Conflict of interest} \\
\small The authors declare no conflict of interest.

\vspace{0.4cm}

\noindent\textbf{\small Data and code availability} \\
\small The code and simulation outputs required to reproduce the analyses are publicly available at \url{https://github.com/hiba-nassor/system-centered-probabilistic-formalism}.

\vspace{0.4cm}

\noindent\textbf{\small Use of AI tools} \\
\small During manuscript preparation, the authors used ChatGPT for language editing, wording suggestions and occasional assistance with code organisation and debugging. All code and scientific content were checked and validated by the authors.

\vspace{0.6cm}
\noindent\rule{\textwidth}{0.4pt}

\section{Introduction}

Ecological and evolutionary systems are often too complex to be characterized by a single set of deterministic parameters alone. Since the seminal work of \citet{May1972}, many studies have addressed this complexity through random interaction matrices, asking probabilistic questions about stability, feasibility and coexistence. This probabilistic viewpoint has become central in the study of complex ecological communities, from random matrix approaches to analyses of high-dimensional community dynamics \citep{AllesinaTang2012, Bunin2017, Barbier2018}.

However, most existing approaches still focus on states one at a time. Analytical studies often ask whether a specific coexistence equilibrium is feasible, saturated, or locally stable under a given interaction structure \citep{May1972, AllesinaTang2012}. In numerical studies, a system is often simulated from one or several initial conditions, and the resulting attractors characterize the model behavior \citep{Strogatz2018, Barbier2018}. Such approaches are essential for identifying dynamically reached outcomes, but they sample the state set through the chosen initial conditions, potentially missing admissible states that are not reached.
This limitation becomes relevant when a system can support several reachable alternative states. Multistability is now recognized as a key element in population dynamics, ecological resilience and regime shifts \citep{Scheffer2001, Beisner2003}. Recent studies further emphasize that alternative stable states are not only theoretical possibilities but can arise experimentally, with consequences for community structure and ecosystem function \citep{unknown, article}. Yet, when equilibria obtained from different systems are pooled and analyzed collectively, the identity of the system that generated each state is lost. As a result, pooling can describe the distribution of states, but not whether a typical system admits no state of a certain type, a unique one, or several alternative states.

Biodiversity, meanwhile, is often treated as a separate question. In empirical ecology, it is quantified through species richness or diversity indices \citep{Tilman1999, Magurran2004}, while in theoretical models it is often reduced to the size of a single feasible or stable community.
What is still missing is a framework that jointly addresses multiplicity, defined as the number of admissible equilibria of a certain type that a system can support, and biodiversity, defined as the number of species coexisting within these equilibria.

Here, we propose a unified probabilistic framework that bridges this gap. Rather than focusing on a single equilibrium or pooled equilibria, we characterize the entire set of admissible states supported by a system.
For a fixed parameter set, we consider several equilibrium types and quantify both their multiplicity and the number of species they contain. We then lift this deterministic description to a probabilistic level by drawing systems at random from a parameter ensemble, treating multiplicity and biodiversity as random variables.
A central element of our approach is the recognition that existence itself is not guaranteed \citep{May1972, Bunin2017}. Some systems admit no admissible state of the chosen type, while others support one or many. This naturally yields a hurdle or zero-inflated statistical structure \citep{Lambert1992}, capturing absence of equilibria, multiplicity and biodiversity within a single distributional object. This formulation makes it possible to disentangle three distinct questions: whether a given state type exists, how many alternative states of that type a system admits, and how many species coexist within them.

We apply this framework to two classes of population models. First, we analyze generalized Lotka-Volterra systems, which govern absolute abundances and form the backbone of theoretical ecology \citep{Lotka1925, Volterra1926}. Second, we study replicator dynamics, which describe how species frequencies change and arise naturally in evolutionary games and multistrain epidemiology \citep{TaylorJonker1978, HofbaeurSigmund1998, MadecGjini2021}. We consider both unstructured random interactions and structured interaction matrices, exploring how added constraints reshapes equilibrium landscapes.
In both cases, the set of admissible equilibria is finite, making the framework computationally possible.

By combining exhaustive equilibrium characterization with probabilistic sampling over systems, our approach reveals emergent patterns that are invisible to state-centered approaches. We show how control parameters reshape not only average biodiversity but also the multiplicity and robustness of coexistence states.

\section{A system-centered probabilistic formalism}

This study is motivated by a shift in perspective on how coexistence and biodiversity are viewed and presented. The aim is to transition from a state-centered perspective, which pools states and ignores systems devoid of states and their multiplicity, to a system-centered perspective, in which the focus is directed towards the system itself, and each system leads to a full probabilistic entity that captures the targeted information. 

For a given system, the outcome of interest is not only how many admissible states a system can support, but also how many species coexist across all admissible states.

This description gives rise to a two-level statistical representation. The first level concerns the probability that a system admits at least one admissible state of a given type (saturated, stable, or evolutionarily stable). The second level concerns the distribution of the number of coexisting species, or possibly extinction, among these states. Indeed, some systems admit no state at all, producing an excess of zeros, while others exhibit variable numbers of coexisting species. 

This motivates hurdle and zero-inflated models, which define a joint probabilistic structure separating existence and multiplicity from the number of coexisting species.

\subsection{Deterministic setting}

We begin by introducing a fully deterministic description and then carry the construction over to a probabilistic setting.

Let $N$ be the total number of species in the system being studied. We define a dynamical system by a set of parameters $p \in \mathcal{P}$, where $\mathcal{P}$ is a parameter space. Depending on the model, $p$ may encode a fitness or an interaction matrix, growth rates, or any additional structural parameters.
Each fixed $p$ defines a completely deterministic system that admits a finite set of states of a given type, denoted by $E(p)$.
All subsequent quantities are defined with respect to $p$.

\begin{definition}[State size] \label{def:state_size}
    For $z \in E(p)$, the size of $z$ is defined as the number of coexisting species, i.e., the cardinality of its support:
    \begin{equation}
        \|z\|_{0} := \mid \mathrm{supp}(z) \mid = \mid \{ i \in \{1,\ldots,N\} \mid z_i \neq 0 \} \mid \, .
    \end{equation}
    
\end{definition}

$\|z\|_{0}$ represents the species richness, which is a basic notion of biodiversity associated to the state $z$. Many alternative indices have been proposed to capture additional aspects of biodiversity such as relative abundances or functional differences between species \citep{Hill, Magurran2004, leincobbo,leinster2024entropydiversityaxiomaticapproach}. These indices will not be considered in the present work.

Next, we define the multiplicity of states supported by a system.

\begin{definition}[Multiplicity] \label{def:multiplicity}
    The multiplicity is defined as the number of admissible states, i.e., the cardinality of $E(p)$:
    \begin{equation}
        S_p := |E(p)| \, .
    \end{equation} 
\end{definition}

This concept appears in various forms, such as multistability in theoretical ecology \citep{Scheffer2001, Beisner2003, Guimera2024}, multiple attractors in nonlinear dynamics \citep{Strogatz2018}, multiple ESS or multiple Nash equilibria in evolutionary game theory \citep{Nash1950, HofbaeurSigmund1998}.
Although they appear in different contexts, all these notions are well captured by the mathematical quantity $S_p$.

To refine the multiplicity by state size, we define the number of states supporting exactly $k$ species.

\begin{definition}[Multiplicity by state size]
    For $ k \in \{0,\dots,N\} $, the multiplicity of states of size $k$ is defined as 
    \begin{equation}
        S^k_p := \mid \{ z\in E(p) \mid \|z\|_0 = k\} \mid \, .
    \end{equation}
\end{definition}

These quantities satisfy $S_p = \displaystyle\sum_{k=0}^N S_p^k$. 
Note that $\{z\in E(p) \mid \|z\|_0=0\} \subset \{0\}$.
Consequently, $S_p^0\in\{0,1\}$: it is equal to one if the extinct state belongs to $E(p)$, and zero otherwise. This case can occur in models whose state space contains the zero vector, such as generalized Lotka-Volterra systems.

In the deterministic setting, multiplicity and biodiversity are therefore fully encoded by the vector
$\left( S^0_p \, , \, \ldots \, , \, S^N_p \right)$.
Under this formulation, we introduce the probability vector in the following way.

\begin{definition}[Probability vector]
        
     Let $\mathcal{E} = \{0,1,\dots,N\} \cup \{\varnothing\}$. The probability vector is  
     \[F_p = \left( F_p^\alpha \right)_{\alpha\in\mathcal{E}} \, ,
     \] 
     with components given, for $\alpha\in\mathcal{E}$, by 
     \[
        F_p^\varnothing:=
        \begin{cases}
            1 & \text{if } S_p = 0 \, ,\\
            0 & \text{otherwise}
        \end{cases}  \quad , \quad \text{ and } \quad \forall k \in \{0,\dots,N\} \quad F_p^k:= \begin{cases}
            0  & \text{if } S_p = 0 \, , \\
            \frac{S_p^k}{S_p} & \text{otherwise}
        \end{cases}
    \]
\end{definition}

\begin{figure}[!htbp]
    \centering
    \includegraphics[width=1\linewidth]{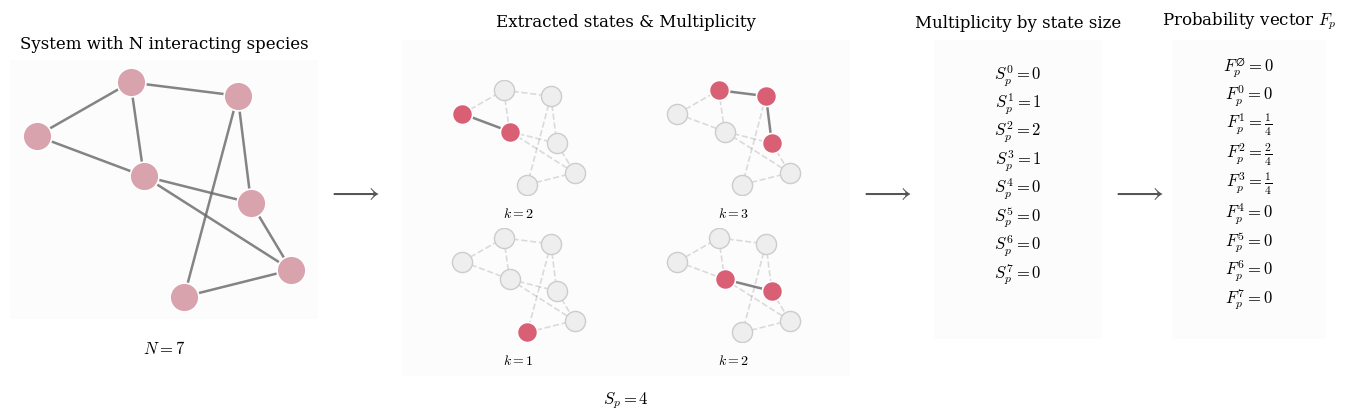}
    \caption{\textbf{Illustration of multiplicity and state size distribution for one parameter set.} Four states of sizes $2$, $3$, $1$, and $2$ are extracted from a system with $N = 7$ interacting species. The total number of states (multiplicity) is $S_p = 4$, and $S_p^k$ counts the number of states of size $k$. These quantities define the probability vector $F_p$, which summarizes the distribution of state sizes for a given system.}
    \label{fig:determ}
\end{figure}

That is, if no state exists ($S_p = 0$), then all probability mass goes to $\varnothing$. Otherwise, $F_p$ normalizes the counts of states of each size into a valid probability distribution over $\mathcal{E}$. Fig.~\ref{fig:determ} illustrates the construction of the probability vector $F_p$ from the states extracted from a single system.

\subsection{Probabilistic formalism} 

With the deterministic structure in place, we now introduce randomness at the level of the parameters. 

Let $(\Omega,\mathcal{F},\mathbb{P})$ be a probability space, the parameter $p$ is no longer fixed; rather, it is viewed as a random variable 
\[
    P : (\Omega, \mathcal{F}, \mathbb{P}) \longrightarrow (\mathcal{P}, \mathcal{B}(\mathcal{P})) \, , 
\]
with distribution $\mathcal{D}$.
All previously defined quantities become random variables when evaluated at $P$.

To capture the resulting randomness in the coexistence outcome, we define the following random variable.

\begin{definition}[Random variable $K^\star$]
    Conditionally on $P=p$, the random variable $K^\star : \Omega \to \mathcal{E}$ \, is defined so that, for all $\alpha \in \mathcal{E}$ ,
    \[
        \mathbb{P}\left(K^\star = \alpha \mid P=p \right) := F_p^{\alpha} \, .
    \]
\end{definition}

In other words, conditional on a parameter realization $p$, $K^\star$ is drawn according to the distribution $F_{p}$.

\begin{definition}[Distribution of multiplicity]
    For $s \in \mathbb{N}$, the distribution of multiplicity $\pi_s$ is defined by 
    \begin{equation}
        \pi_s := \mathbb{P}(S_P = s) \, .
        \label{eq:pi_s}
    \end{equation}
\end{definition}

The distribution $\pi_s$ quantifies how often systems possess no admissible state ($s=0$), a unique one ($s=1$), or multiple ones ($s\ge2$), which leads to a richer understanding of the structural robustness of the systems.

\begin{definition}[Conditional distribution]
    For $k \in \{0,\ldots,N\}$, the conditional distributions of $K^\star$ are defined as follows: 

    Conditionally on $\{S_P > 0\}$,
    \begin{equation}
        p_k :=  \mathbb{P}\left(K^\star = k \mid S_P>0 \right) \, ,
    \label{eq:p_k}
    \end{equation}

    and, for $s \in \mathbb{N^\ast}$, conditionally on $\{S_P = s\}$,
    \begin{equation}
        p_{k \mid s} := \mathbb{P}\left(K^\star = k \mid S_P = s \right) 
    \label{eq:p_ks}
    \end{equation}

\end{definition}

\begin{remark}
    These distributions can be expressed in terms of the probability vector $F_P$. In particular,
    \[
        p_k = \mathbb{E}[F_P^k \mid S_P > 0], \quad p_{k \mid s} = \mathbb{E}[F_P^k \mid S_P = s].
    \]
\end{remark}

Based on these considerations, the random variable $K^\star$ naturally decomposes as

\begin{equation}
    \boxed{
    \mathbb{P}(K^\star = \alpha) =
    \begin{cases}
    \pi_0, & \text{if } \alpha = \varnothing,\\
    (1 - \pi_0)\, p_k, & \text{if } \alpha = k,\ \ k \in \{0,\ldots,N\}
    \end{cases}
    }
    \label{eq:hurdle_model}
\end{equation}

This is the explicit form of a hurdle model \citep{MULLAHY1986341}: a Bernoulli variable decides whether a state exists, and conditionally on existence, the system follows a distribution over the number of coexisting species.

Conditionally on existence, the expected number of coexisting species is
\begin{equation}
    \mathbb{E}\left[K^\star \mid S_P>0\right] = \sum_{k=0}^{N} k\,p_k \, .
\label{eq:expected_state_size}
\end{equation}
This quantity measures the mean size of an existing state, after first conditioning on state existence and then weighting states uniformly within each system.

Rather than conditioning on $\{S_P > 0\}$, we may instead condition on $S_P = s$, which is equivalent to focusing on the systems that admit exactly $s$ states. Hence, the law of $K^\star$ can then be decomposed over the multiplicity levels as

\begin{equation}
    \boxed{
    \mathbb{P}(K^\star = \alpha)
    =
    \begin{cases}
    \pi_0 ,
    & \text{if } \alpha = \varnothing,
    \\
     \sum_{s\ge1} \pi_s\, p_{k \mid s},
    & \text{if } \alpha = k,\ \ k \in \{0,\ldots,N\}
    \end{cases}
    }
    \label{eq:decom_mult}
\end{equation}

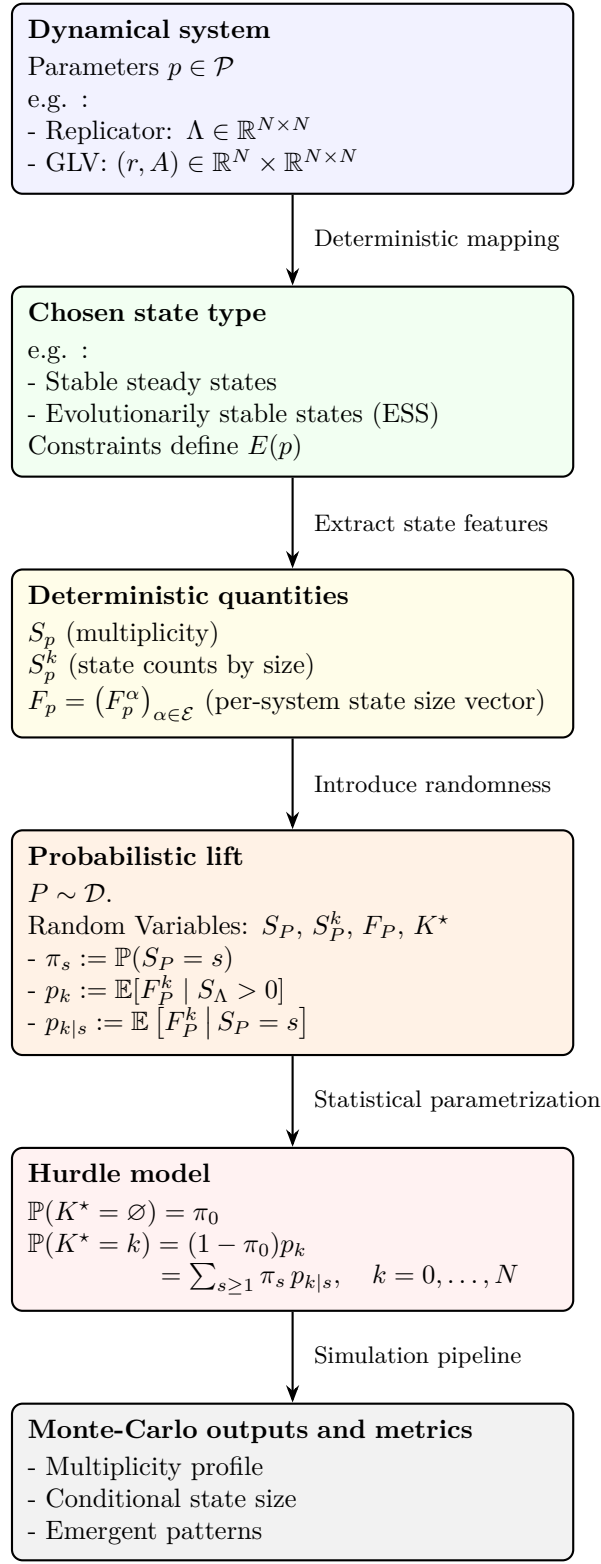
\begin{figure}[!htbp]
\centering

\begin{tikzpicture}[
    font=\small,
    block/.style={
        rectangle,
        rounded corners,
        draw=black,
        thick,
        align=left,
        inner sep=6pt,
        text width=7cm
    },
    blocksys/.style={block, fill=blue!5},
    blockeq/.style={block, fill=green!5},
    blockdet/.style={block, fill=yellow!10},
    blockprob/.style={block, fill=orange!10},
    blockmix/.style={block, fill=red!5},
    blockout/.style={block, fill=gray!10},
    arrow/.style={->,>=Stealth,thick}
]

\matrix (m) [row sep=1.2cm] {

    \node[blocksys] (system) {%
        \textbf{Dynamical system}\\[2pt]
        Parameters $p\in\mathcal{P}$\\
        e.g. :\\
        - Replicator: $\Lambda\in\mathbb{R}^{N\times N}$\\
        - GLV: $(r,A)\in\mathbb{R}^{N}\times\mathbb{R}^{N\times N}$
    }; \\

    \node[blockeq] (eqtype) {%
        \textbf{Chosen state type}\\[2pt]
        e.g. : \\
        - Stable steady states \\ 
        - Evolutionarily stable states (ESS)\\
        Constraints define $E(p)$
    }; \\

    \node[blockdet] (det) {%
        \textbf{Deterministic quantities}\\[2pt]
        $S_p $ (multiplicity) \\
        $S_p^k$ (state counts by size)\\
        $F_p = \left( F_p^\alpha \right)_{\alpha\in\mathcal{E}}$ (per-system state size vector)
    }; \\

    \node[blockprob] (prob) {%
        \textbf{Probabilistic lift}\\[2pt]
         $P\sim\mathcal{D}$.\\
        Random Variables: $S_P$, $S_P^k$, $F_P$, $K^\star$\\
        - $\pi_s := \mathbb{P}(S_P=s)$\\
        - $p_k := \mathbb{E}[F_P^k\mid S_\Lambda>0]$\\
        - $p_{k \mid s} := \mathbb{E}\left[ F_P^k \,\big|\, S_P = s \right]$
    }; \\

    \node[blockmix] (mixture) {%
        \textbf{Hurdle model}\\[2pt]
        $\mathbb{P}(K^\star=\varnothing)=\pi_0$ \\
        $
        \mathbb{P}(K^\star=k) = (1 - \pi_0) p_k$\\
        $\qquad \qquad \quad = \sum_{s\ge1} \pi_s \, p_{k \mid s}, \quad k = 0,\ldots,N
        $
    }; \\

    \node[blockout] (output) {%
        \textbf{Monte-Carlo outputs and metrics}\\[2pt]
        - Multiplicity profile\\
        - Conditional state size\\
        - Emergent patterns
    }; \\
};

\draw[arrow] (system) -- node[right=4pt]{\footnotesize Deterministic mapping} (eqtype);
\draw[arrow] (eqtype) -- node[right=4pt]{\footnotesize Extract state features} (det);
\draw[arrow] (det) -- node[right=4pt]{\footnotesize Introduce randomness} (prob);
\draw[arrow] (prob) -- node[right=4pt]{\footnotesize Statistical parametrization} (mixture);
\draw[arrow] (mixture) -- node[right=4pt]{\footnotesize Simulation pipeline} (output);

\end{tikzpicture}

\caption{\textbf{Conceptual workflow summarizing the probabilistic formalism}. A deterministic dynamical system is specified by parameters $p$, from which a chosen state type defines the admissible set $E(p)$ and the associated deterministic quantities $(S_p,S_p^k,F_p)$. Randomizing the parameters ($P\sim\mathcal{D}$) lifts these objects to random variables and induces a hurdle or zero-inflated model for the state size $K^\star$, from which presence or absence of coexistence, multiplicity, and emergent coexistence patterns can be quantified.}
\label{fig:workflow}
\end{figure}

This representation sheds light on how the probability of observing a state of size $k$ results from two components: the probability $\pi_s$ that a system admits exactly $s$ state of the prescribed type, and the internal composition $p_{k \mid s}$ of such system.
Fig.~\ref{fig:workflow} summarizes the transition from the deterministic characterization of a dynamical system to its probabilistic representation.

\begin{remark}
    In models such as the generalized Lotka-Volterra system, the fully extinct state $z\equiv0$ can itself be a valid state.
    By merging the absence of states ($\varnothing$) and the extinction state ($k=0$) into a single case, we obtain a variable denoted $K\in\{0,1,\dots,N\}$ with law 
    
    \begin{equation}
        \mathbb{P}(K = k) =
        \begin{cases}
        \pi_0 + (1 - \pi_0)\, p_0, & \text{if } k = 0,\\
        (1 - \pi_0)\, p_k, & \text{if } k \in \{1,\ldots,N\}
        \end{cases}
        \label{eq:zero_inflated_model}
    \end{equation}
    
    This compressed version corresponds to \textbf{zero-inflated model} \citep{Lambert1992, Zuur2009, WengerFreeman2008}, where zeros arise from two different sources: a structural zero corresponding to the absence of a state, and a dynamical zero corresponding to an extinction state. It also allows standard numerical summaries, such as the mean and variance of state size, to be defined without conditioning on the existence of a state.
    
\end{remark}

\subsection{Comparison with the state-centered approach}

We now make explicit the distinction between the distribution of state sizes from a state-centered view, denoted by $p_k^{\mathrm{state}}$, and the system-centered view, denoted by $p_k^{\mathrm{system}}$.
Both descriptions are considered conditionally on the existence of admissible states.
Using the same notation as before, we recall that
\begin{equation}
    p_k^{\mathrm{system}} = \mathbb{E}[F_P^k \mid S_P>0] = \mathbb{E}\left[ \frac{S_P^k}{S_P} \mid S_P > 0 \right] \, .
\end{equation}  
We now consider the other case where, given $m$ i.i.d. systems $\left(P^{(i)}\right)_{ 1\leq i \leq m }$, the distribution $\widehat{p}^{\mathrm{state}} = \left(\widehat{p}^{\mathrm{state}}_k\right)_k$ is defined by 
\begin{equation*}
    \widehat{p}_k^{\mathrm{state}} = \frac{\displaystyle\sum_{i : S_{P^{(i)}} > 0} S_{P^{(i)}}^k}{\displaystyle\sum_{i : S_{P^{(i)}} > 0} S_{P^{(i)}}} \, .
\end{equation*}
By the law of large numbers applied to systems satisfying $S_P > 0$, as $m \rightarrow \infty$ this converges almost surely to
\begin{equation}
    p_k^{\mathrm{state}} = \frac{\mathbb{E}[S_P^k \mid S_P > 0]}{\mathbb{E}[S_P \mid S_P > 0]} \, .
\end{equation}
Thus, this distribution corresponds to sampling a state uniformly at random among all states generated by all systems. Using the identity ($S_p^k = S_p F_p^k$) and expanding the expectation of this product according to the definition of covariance, we obtain the following
\begin{align*}
    p_k^{\mathrm{state}} & = \frac{\mathbb{E}[S_P \mid S_P > 0 ] \cdot \mathbb{E}[F_P^k \mid S_P > 0] + \mathrm{Cov}\left(S_P, F_P^k \mid S_P > 0\right)}{\mathbb{E}[S_P \mid S_P > 0]} \, .
\end{align*}
Linking the two approaches yields the following result
\begin{equation}
    \boxed{
    p_k^{\mathrm{state}} \, = \, p_k^{\mathrm{system}}  \, +  \, \frac{\mathrm{Cov}\left(S_P, F_P^k \mid S_P > 0\right)}{\mathbb{E}[S_P \mid S_P > 0]}
    }
    \label{eq:cov}
\end{equation}
Thus, the difference between the two distributions is entirely governed by the covariance term, which captures the interaction between multiplicity (through $S_P$) and the internal structure of states (through $F_P^k$).

In the context of multistable dynamical systems, the quantity $S_P$ is itself a central object of interest. The state-centered approach mixes two distinct effects: the number of states per system and the composition of these states.
As a consequence, systems with a large number of states dominate the statistics, potentially masking the behavior of typical systems.
In contrast, the system-centered approach preserves the system-level variability and allows one to disentangle the existence of states, their multiplicity, and their composition. It therefore provides a more appropriate framework for analyzing how multistability depends on model parameters.

\subsection{Two case studies} 
To illustrate the applicability of the formalism to population dynamics models, we consider two classical families of interacting species systems: the generalized Lotka-Volterra (GLV) system and the replicator dynamics. These two models are closely related: an \(n\)-strategy replicator equation can be rewritten as a GLV system with \(n-1\) variables, and conversely GLV dynamics can be embedded in a replicator framework after a suitable change of variables. This dynamical connection motivates treating these two models, as they provide distinct biological interpretations, but closely related mathematical structures \citep{HofbaeurSigmund1998, Allesina2026}.

In both models, the parameter space $\mathcal{P}$ can be associated to a $\sigma$-algebra, forming a measurable space.
Each choice of a subcommunity leads to at most one equilibrium satisfying the corresponding linear constraints. Since the number of possible subcommunities is finite, the equilibrium set $E(p)$ is also finite.  

This guarantees that the vector $(S_p^k)_{1\leq k\leq N}$ and the random variable $K^\star$ are well-defined, and the probabilistic formalism developed above applies directly.

We search for states satisfying admissibility, saturation, and asymptotic stability. Moving beyond classical stability, we also explore the notion of evolutionary stability in the context of replicator dynamics.

\section{Generalized Lotka-Volterra systems}

\subsection{Model definition and equilibrium types}
We consider the generalized Lotka-Volterra (GLV) model, a classical framework for the dynamics of multispecies population \citep{HofbaeurSigmund1998, Takeuchi}.
The GLV model describes the evolution of absolute abundances $x =(x_1,\dots,x_N) \in\mathbb{R}_+^N$ with intrinsic growth vector $r \in \mathbb{R}^N$ and interaction matrix $A \in \mathbb{R}^{N \times N}$:
\begin{equation}
    \dot x_i = x_i\big(r_i + (Ax)_i\big), \quad \forall i \in \{1, \ldots, N\}.
    \label{eq:GLV}
\end{equation}

\begin{definition}[Admissible steady state]\label{def:s_glv}
    $x^\ast \in \mathbb{R}_+^N$ is an \emph{ admissible steady state} if it satisfies
    \begin{equation}
        x_i^\ast \bigl( r_i + (A x^\ast)_i \bigr) = 0, \quad \forall i \in \{1, \ldots, N\}.
        \label{eq:s_glv}
    \end{equation}

\end{definition}

Generally, a steady state $x^\ast$ is given by setting \eqref{eq:GLV} to zero, i.e. by either $x^\ast_i = 0$ for extinct species, or $r_i + (Ax^\ast)_i = 0$ for surviving ones, which results generically in $2^N$ possible outcomes. To be admissible, this equilibrium has to satisfy $x_i^\ast \geq 0$ for all $i \in \{1, \ldots, N\}$. However, the extinction equilibrium $x^\ast=0$ is always feasible, which introduces an additional equilibrium of size $0$.

\begin{definition}[Saturated state]\label{def:sat_glv}
    A steady state $x^\ast \in \mathbb{R}_+^N$ is said to be \emph{saturated} if it satisfies
    \begin{equation}
        r_i + (A x^\ast)_i \le 0, \quad \forall i \in \{1, \ldots, N\}.
        \label{eq:sat_glv}
    \end{equation}
\end{definition}

Among all admissible steady states, saturated states are those that cannot be invaded by any absent species, i.e $x_i^\ast = 0 $. The condition \eqref{eq:sat_glv} ensures that the growth rate of any absent species is non-positive \citep{HofbaeurSigmund1998}.

\begin{definition}[Asymptotically stable state]\label{def:sta_glv}
    A steady state $x^\ast \in \mathbb{R}_+^N$ is said to be \emph{asymptotically stable} if all the eigenvalues of its Jacobian matrix have negative real parts.
\end{definition}

The Jacobian matrix of \eqref{eq:GLV} evaluated in steady state has a block-triangular structure. 
Let $I \subset \{1,\dots,N\}$ denote the support of the steady state, with $x_i^\ast>0$ for $i\in I$ and $x_i^\ast=0$ for $i\in I^c$. 
After reordering the variables so that the species in $I$ appear first, the Jacobian can be written as

\begin{equation}
    J(x^\ast)
    =
    \begin{pmatrix}
    \operatorname{diag}(x_I^\ast)A_{II}
    &
    \operatorname{diag}(x_I^\ast)A_{I I^c}
    \\
    0
    &
    \operatorname{diag}\!\left(r_{I^c}+A_{I^c I}x_I^\ast\right)
    \end{pmatrix}.
    \label{eq:glv_jacobian_block}
\end{equation}
One block governs the internal dynamics of species present at the state, while the other block is diagonal and corresponds to the growth rates of absent species. Since a block-triangular matrix has a spectrum equal to the union of the spectra of its diagonal blocks \citep{Horn_Johnson_1985}, asymptotic stability requires that all invasion growth rates be strictly negative. Thus the saturation condition
\eqref{eq:sat_glv} is necessary for stability, and can be used to exclude many unstable steady states without computing the full Jacobian spectrum.

\subsection{Interaction structures and numerical procedure}

Although the equilibrium types described above arise independently of the choice of parameters, their distribution and prevalence depend strongly on the structure of the parameters. To investigate how different interaction shape the probabilistic outputs, we consider these two particular configurations.

\paragraph{Random competitive. }In the first interaction structure, growth rates $r_i$ are drawn independently from a standard normal distribution and are constrained to be positive. Interaction coefficients $A_{ij}$ are drawn independently from a centered Gaussian distribution and constrained to be non-positive, corresponding to purely competitive interactions.

\paragraph{Structured $(c,\sigma)$. }In the second interaction structure, the growth rates are fixed to $r_i=1$, and the interaction matrix is defined as
\begin{equation*}
    A = -I_N + B,
    \label{eq:structure_A}
\end{equation*}
where $I_N$ is the identity matrix and $B$ is a random matrix with zero diagonal. For $i\neq j$, the entries of $B$ are given by
\begin{equation*}
    B_{ij} = X_{ij} U_{ij},
\end{equation*}
where
\[
X_{ij} \sim \mathrm{Bernoulli}(c),
\qquad
U_{ij} = -|Z_{ij}|,\quad Z_{ij}\sim \mathcal{N}(0,\sigma^2).
\]
The parameter $c\in[0,1]$ controls the connectance of the interaction network, while $\sigma>0$ controls the typical interaction strength. 

These two structures are standard variants of random Lotka-Volterra community models. The first corresponds to randomly assembled competitive GLV systems, while the
second is a sparse competitive May-type matrix with negative self-regulation, connectance $c$, and interaction strength $\sigma$
\citep{May1972,GOH197763,AllesinaTang2012,Bunin2017}.

\begin{table}[H]
\centering
\caption{\textbf{Summary of the interaction structures considered in the GLV model.} }
\begin{tabular}{|l|l|l|}
\hline
Structure & Growth rates $r_i$ & Interaction matrix $A$ \\
\hline
Random competitive
& $r_i = |Y_i|, \; Y_i\sim \mathcal{N}(0,1)$
& $A_{ij} = -|G_{ij}|, \; G_{ij}\sim \mathcal{N}(0,1)$ \\[4pt]
\hline

Structured $(c,\sigma)$
& $r_i = 1$ 
& $A=-I_N+B,\; \text{with } \,B_{ii}=0,$ \\
& & $\text{for } i\neq j , \; B_{ij}=X_{ij}U_{ij},\; X_{ij}\sim\mathrm{Bernoulli}(c),$ \\
& & $U_{ij}=-|Z_{ij}|,\; Z_{ij}\sim\mathcal{N}(0,\sigma^2)$ \\

\hline
\end{tabular}
\label{tab:GLV_structures}
\end{table}

\subsection{Results}

We now illustrate the output of the formalism for the competitive GLV model. 
For each interaction structure described in Table~\ref{tab:GLV_structures}, we consider systems of $N=10$ interacting species and generate  $M=10\,000$ independent systems. 
For each sampled parameter set $p$, all $2^N$ possible supports are enumerated, and the corresponding admissible states are classified as saturated or stable according to the criteria introduced above. 
This produces, for each system $p$, a set of admissible states of the chosen type $E(p)$.

From this set, we extract three complementary pieces of information: the multiplicity $S_p$, defined as the number of admissible states of a certain type associated with system $p$, the state size distribution within each system, summarized by the probability vector $F_p$, and the conditional distribution of state sizes given the multiplicity $ p_{k \mid s} = \mathbb{P}(K^\ast=k\mid S_P=s)$, which describes how state sizes are distributed among systems with the same multiplicity class.

These quantities are the basis of the statistical summaries used in the following: the multiplicity distribution $(\pi_s)_s$ in Eq.~\eqref{eq:pi_s}, the state size distributions for each multiplicity class $(p_{k\mid s})_k$ in  Eq.~\eqref{eq:p_ks}, the distribution of $K^\star$ in  Eq.~\eqref{eq:hurdle_model}, and the conditional mean state size $\mathbb{E}(K^\star\mid S_P>0)$ in Eq.~\eqref{eq:expected_state_size}.

The following results are organized according to the type of state considered. 
We first analyze saturated states, we then restrict the analysis to stable states. In both cases, the figures should be read as different projections of the same formalism: existence, multiplicity, and state size distribution.

\subsubsection{Saturated states}
\label{sec:glv_saturated_states}
 
Among admissible states, saturated states describe the property of noninvadability (Def.~\ref{def:sat_glv}), before imposing internal dynamical stability.

\begin{figure}[!htbp]
    \centering
    \includegraphics[width=\textwidth]{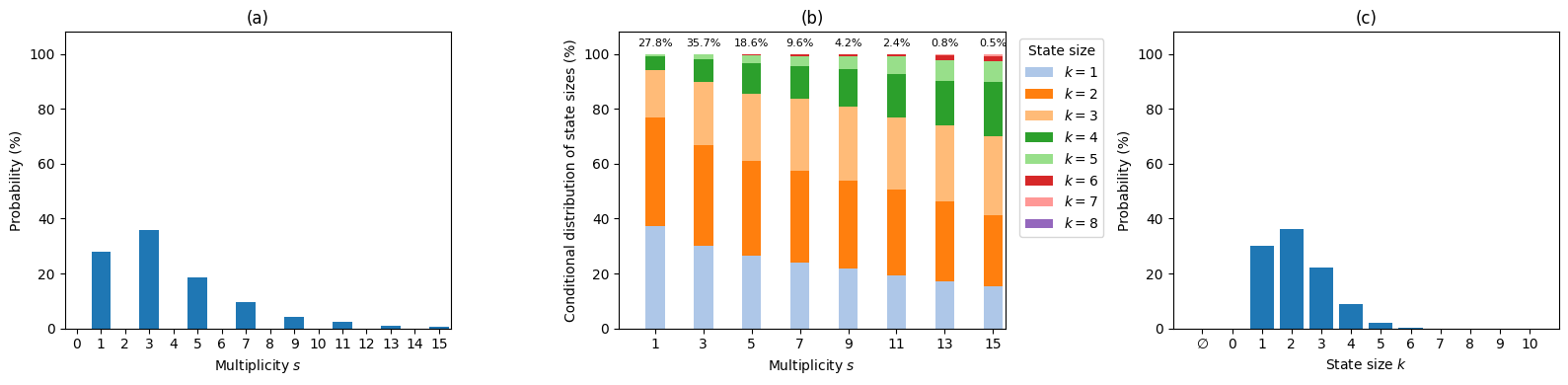}
    \caption{
    \textbf{Output of the formalism for saturated states in the random competitive GLV model.} The random competitive structure is defined in Table~\ref{tab:GLV_structures}.
    Saturated states are identified after enumerating all possible supports.
    (a) Distribution of the number of saturated states $s$ (multiplicity).
    (b) Conditional distribution of state sizes for each multiplicity value $s$. Each bar is normalized to $100\%$, the percentages displayed above the bars indicate the fraction of systems belonging to that multiplicity class, matching the distribution shown in Pannel~(a). Each color corresponds to a state size observed at least once, showing that each multiplicity class contains saturated states of many different sizes.
    For clarity, Panels~(a)--(b) display only multiplicity classes up to $s_{\max}=15$. 
    (c) Distribution of state sizes $k$ obtained by aggregating over all systems and multiplicities.
    Together, Panels~(a)--(c) provide marginal summaries of the information jointly combined in Panel~(b).}
    \label{fig:glv_sat_random_formalism}
\end{figure}

Fig.~\ref{fig:glv_sat_random_formalism} gives a first overview of the distribution of saturated states in the random competitive ensemble. 
Across systems, multiplicity varies significantly, with the existence always guaranteed. In particular, only odd multiplicity values are observed. Let us emphasize that theoretical results from linear complementarity problem prove the existence of saturated states and the oddness of their multiplicity \citep{MURTY197265, clenet}. 
The conditional distribution of state sizes differs: low multiplicity tends to host states with smaller number of coexisting species, whereas higher multiplicities are associated with the appearance of large state sizes; though this trend should be interpreted with caution, as high multiplicity systems represent only a small fraction of the ensemble. At the aggregate level, the distribution of the number of coexisting species peaks at $k = 2$, and drops sharply beyond $k = 4$.

\begin{figure}[!htbp]
    \centering
    \includegraphics[width=\textwidth]{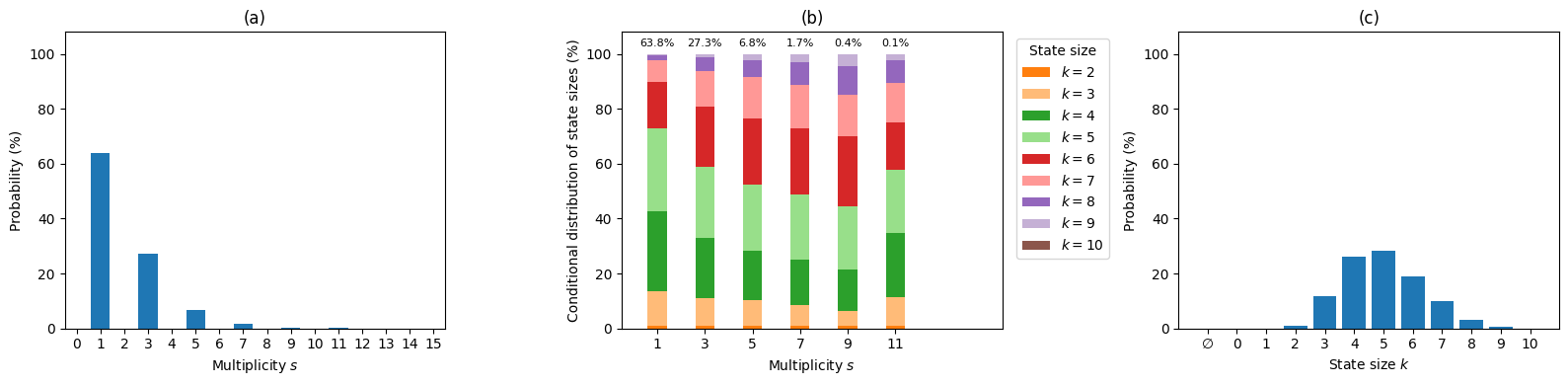}
    \caption{
    \textbf{Output of the formalism for saturated states in the structured GLV model for fixed connectance and interaction strength.} The structure of the interaction parameters is defined in Table~\ref{tab:GLV_structures}, with $c=0.9$ and $\sigma=0.8$.
    Saturated states are identified after enumerating all possible supports.
    The panels follow the same conventions as in Fig.~\ref{fig:glv_sat_random_formalism}: Panel~(a) shows the distribution of the number $s$ of saturated states, Panel~(b) shows the conditional distribution of state sizes within each multiplicity class, and Panel~(c) shows the marginal distribution of state sizes obtained by aggregating over all saturated states. 
    }
    \label{fig:glv_sat_structured_formalism}
\end{figure}

Fig.~\ref{fig:glv_sat_structured_formalism} defines a competitive structure, with $c=0.9$ and $\sigma=0.8$. Compared to the totally random competitive ensemble, this regime shifts the outcome towards low multiplicity: systems with a single saturated state now dominate, accounting for $64\%$ of cases, while saturated states themselves tend to host more species, with the distribution peaking around $k = 5$.
Panel~(b) reveals that the conditional size distribution remains highly diverse across all multiplicity classes, with no states with one species. This suggests that interaction strength and connectance shape saturation outcomes in a non-trivial way; not only by shifting the typical number of saturated states, but also by reorganizing their internal composition.

\begin{figure}[!htbp]
    \centering
    \includegraphics[width=\textwidth]{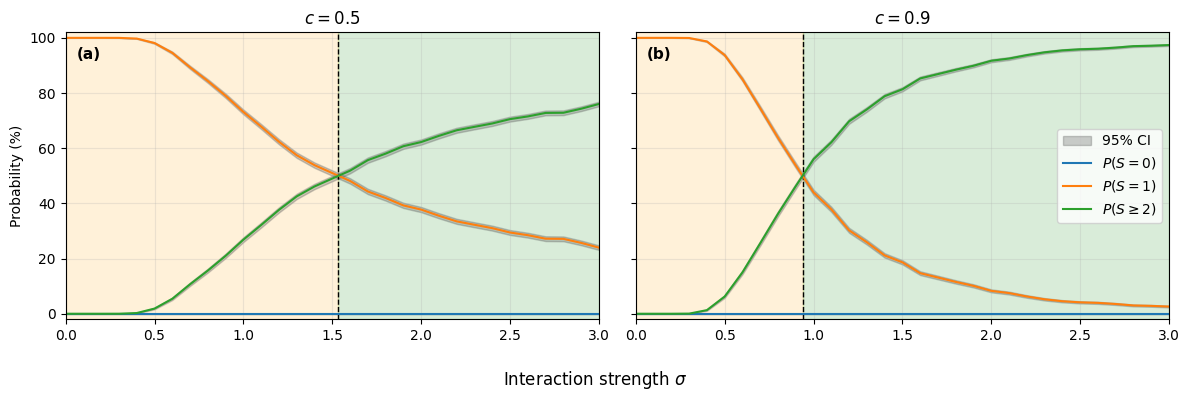}
    \caption{
    \textbf{Multiplicity regimes of saturated states as a function of interaction strength in the structured GLV model.}
    For each value of the interaction strength $\sigma \in [0,3]$, with step size $0.1$, systems are generated using the structured interaction ensemble defined in Table~\ref{tab:GLV_structures}, for two fixed connectance values: (a) $c=0.5$ and (b) $c=0.9$.
    The curves show the probabilities of observing no saturated state, $P(S=0)$, exactly one saturated state, $P(S=1)$, or multiple saturated states, $P(S\geq2)$.
    The grey bands indicate approximate $95\%$ confidence intervals for the Monte-Carlo estimates.
    The shaded background regions indicate the dominant multiplicity regime: orange where $P(S=1)$ is largest, and green where $P(S\geq2)$ is largest.
    In each panel, the vertical dashed line marks the transition point at which $P(S=1)$ and $P(S\geq2)$ are equal.
    }
    \label{fig:glv_sat_mono_multi}
\end{figure}

We then isolate the multiplicity component by tracking the probabilities of absence, uniqueness, and multiplicity as interaction strength $\sigma$ and connectance $c$ vary. We illustrate the cases of $c=0.5$ and $c=0.9$ in Fig.~\ref{fig:glv_sat_mono_multi}

Fig.~\ref{fig:glv_sat_mono_multi} shows how the multiplicity of saturated states changes with the interaction strength $\sigma$ for two fixed connectance values. In both cases, the probability of observing no saturated state remains equal to zero over the whole range of $\sigma$. The main transition is therefore not between the absence and presence of saturated states, but between uniqueness and multiplicity. For weak interactions ($\sigma$ small), almost all systems admit a single saturated state. As $\sigma$ increases, this regime is replaced by a regime characterized by several saturated states. This transition is marked by the intersection between $P(S=1)$ and $P(S\geq 2)$. Increasing connectance shifts this transition to smaller interaction strengths. For $c=0.5$, the uniqueness of saturated states persists over a broader range of $\sigma$, and the crossover occurs only around the intermediate interaction strengths. In contrast, for $c=0.9$, multiplicity emerges much earlier: the probability of multiple saturated states rapidly increases and dominates for most of the explored range.

\begin{figure}[!htbp]
    \centering
    \includegraphics[width=0.75\linewidth]{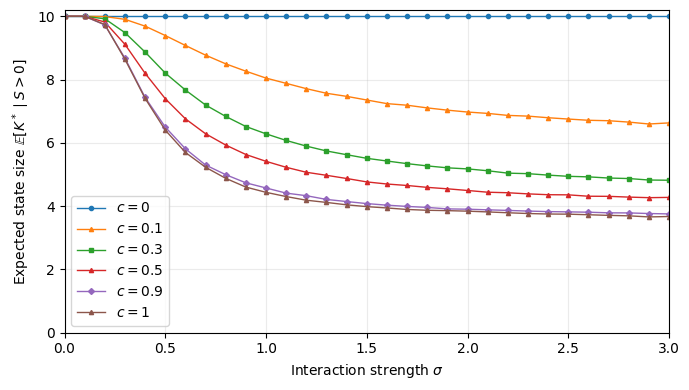}
    \caption{
    \textbf{Expected number of coexisting species conditionally on the existence of saturated states as a function of interaction strength $\sigma$ in the structured GLV model for fixed connectance values.} For each value of the interaction strength $\sigma \in [0,3]$, with step size $0.1$, independent systems are generated using the structured interaction ensemble defined in Table~\ref{tab:GLV_structures}, for several fixed connectance values: $c=0$, $c=0.1$, $c=0.3$, $c=0.5$, $c=0.9$ and $c=1$. The curves show the conditional expectation $\mathbb{E}[K^\ast \mid S>0]$, where $K^\ast$ is the size of a saturated state sampled within a system and $S$ is the number of saturated states in that system. 
    }
    \label{fig:glv_sat_exp}
\end{figure}

Multiplicity tells us how many saturated states a system admits but not their composition. This information is given by the conditional mean state size $\mathbb{E}[K^* \mid S > 0]$ (Fig.~\ref{fig:glv_sat_exp}), computed among systems admitting at least one saturated state.

In the disconnected case, $c=0$, the conditional expected size remains equal to the full system size, $\mathbb{E}[K^\ast\mid S>0] = N$, over the whole range of $\sigma$. This provides a useful reference case: when the species do not interact, the interaction matrix $A$ is equal to the negative identity matrix $-I_N$, independently of $\sigma$. Thus, the dynamics are reduced to $N$ independent logistic equations and the only saturated state is the one containing all species independently of $\sigma$. The same reasoning applies at $\sigma=0$: $A = - I_N$ independently from $c$, and all curves start from the maximal value $\mathbb{E}[K^\ast\mid S>0]=N$. 

As soon as connectance becomes positive, increasing $\sigma$ leads to a monotonic decrease in the mean size of saturated states, indicating that stronger interactions reduce the typical number of coexisting species within saturated states.

For low connectance, the expected size declines gradually and remains relatively large even at high interaction strength. In contrast, for high connectance, the decline is much sharper, and the curves approach lower plateau values. Together with Fig.~\ref{fig:glv_sat_mono_multi}, this shows that increasing $\sigma$ and $c$ has a dual effect: it promotes the emergence of multiple saturated states while reducing the number of coexisting species.

\subsubsection{Stable states}
\label{sec:glv_stable_states}

Stability adds a second filter to the saturated property. 
A stable state must be non-invadable by absent species, but it must also be internally stable with respect to small perturbations of the species present in its support. 

\begin{figure}[!htbp]
    \centering
    \includegraphics[width=\textwidth]{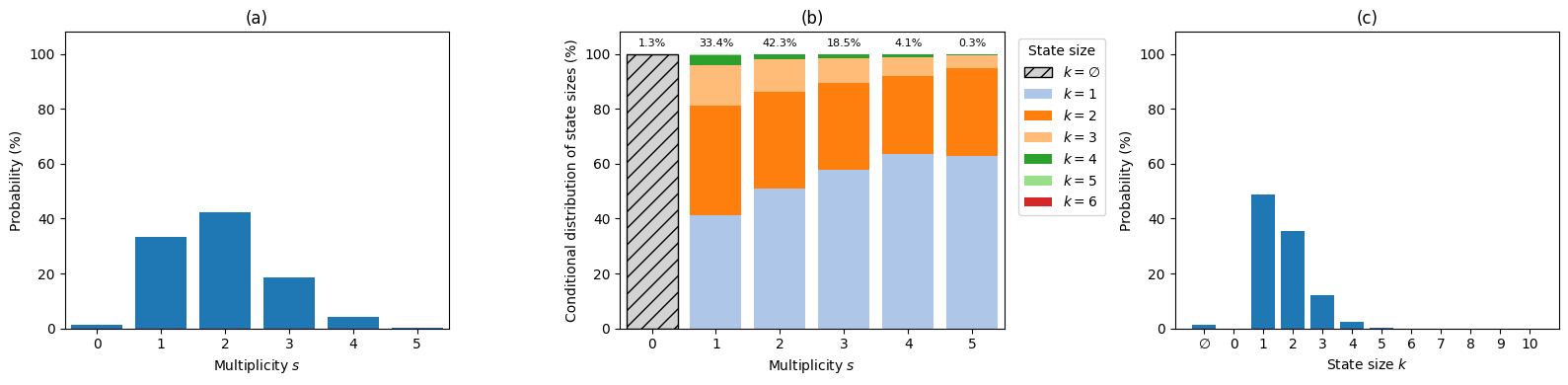}
    \caption{
    \textbf{Output of the formalism for stable states in the random competitive GLV model.} The panels show the multiplicity distribution, the conditional distribution of state sizes by multiplicity class, and the marginal distribution of state sizes distribution. The hatched bar corresponds to systems with no stable state. Same sampling procedure and graphical conventions as in Fig.~\ref{fig:glv_sat_random_formalism}, but applied to stable states. 
    }
    \label{fig:glv_sta_random_formalism}
\end{figure}

\begin{figure}[!htbp]
    \centering
    \includegraphics[width=\textwidth]{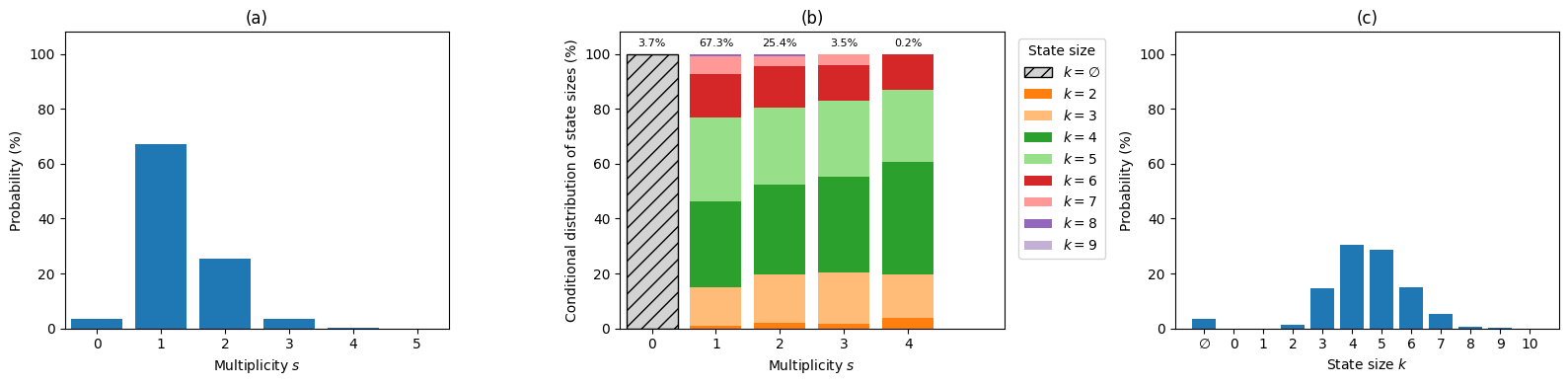}
    \caption{
    \textbf{Output of the formalism for stable states in the structured GLV model for fixed connectance and interaction strength.}  The Panels show the multiplicity distribution, the conditional distribution of state sizes by multiplicity class, and the marginal distribution of state sizes distribution. Same parameter values, sampling procedure and graphical conventions as in Fig.~\ref{fig:glv_sat_structured_formalism}, with $c=0.9$ and $\sigma=0.8$, but applied to stable states.
    }
    \label{fig:glv_sta_structured_formalism}
\end{figure}

This is immediately visible in Fig. ~\ref{fig:glv_sta_random_formalism}--~\ref{fig:glv_sta_structured_formalism}. Compared to the saturated state case, the multiplicity distribution is compressed toward smaller values of $s$: systems with many saturated states often retain only few stable states. Additionally, a fraction of systems admit no stable state at all; the $s=0$ class, absent in saturated states, now accounts for a visible proportion of systems.

Among systems that admit stable states, the multiplicity remains low: monostability ($s=1$) and bistability ($s=2$) together dominate, with higher multiplicities becoming progressively rarer. Compared to saturated states, the distribution is therefore more concentrated at low multiplicities, and the oddness constraint no longer applies. Panel~(b) reveals a contrasting pattern with respect to the saturated case: as $s$ increases, stable states tend to involve fewer coexisting species, with $k=1$ becoming increasingly dominant. The state size distribution is strongly shifted toward a small number of coexisting species, with $k=1$ and $k=2$ accounting for the majority of stable states, a direct consequence of the stability filter, which tends to eliminate larger, more fragile equilibria. Therefore, saturation and stability act as two different filters. Saturation is an external condition: the species outside the support must have a negative growth rate. In a competitive GLV system with a fixed number of species, larger state sizes leave fewer absent species to exclude. This makes larger states more likely to satisfy the saturation condition. Stability, however, requires an additional internal condition: the species present must remain locally stable within the state support, excluding more complex dynamics. When the support size increases, there are more interactions to balance, so larger states are more easily destabilized. Thus, a saturated state may exclude absent species, but it is stable only if the present species can coexist in a stable way.

Fig.~\ref{fig:glv_sta_structured_formalism} applies the same structured regime as before ($c=0.9$, $\sigma=0.8$), now through the lens of stable states. The stability filter considerably reshapes the results relative to the saturated case. Monostability now dominates; $s=1$ accounts for $67.3\%$ of systems, while the fraction of systems admitting no stable state at all rises noticeably, reaching $3.7\%$. The overall multiplicity is reduced, comparing to the high values of $s$ reached in the saturated case. Panel~(b) shows that, conditional on multiplicity, the stable states tend to be intermediate to large size, with $k=4$, $5$ well represented across multiplicity classes, strikingly different from the random case, where the stable states were concentrated at $k=1$. This is confirmed by Panel~(c), where the size distribution peaks around $k=4$, reflecting the dense competitive regime in which larger coexisting species are found. The contrast between these two configurations is not only quantitative, but also reflects a change in the balance between interspecific competition and intraspecific regulation. In the random competitive case, diagonal and off-diagonal entries are sampled on the same scale, so self-regulation is not always stronger than competition from other species. Therefore, coexistence becomes more fragile as the state size grows, and the stability constraints mainly preserve small stable states. On the contrary, in the structured case the decomposition $A = -I_N + B$ imposes a fixed intraspecific regulation, while interspecific interactions are controlled by the connectance $c$ and the interaction strength $\sigma$. With these parameter values, stronger self-regulation makes coexistence easier, so stable states can contain more species. This interpretation is consistent with classical complexity-stability arguments and with GLV coexistence theory, where coexistence is favored when intraspecific competition dominates interspecific competition \citep{May1972, AllesinaTang2012}.

\begin{figure}[!htbp]
    \centering
    \includegraphics[width=\textwidth]{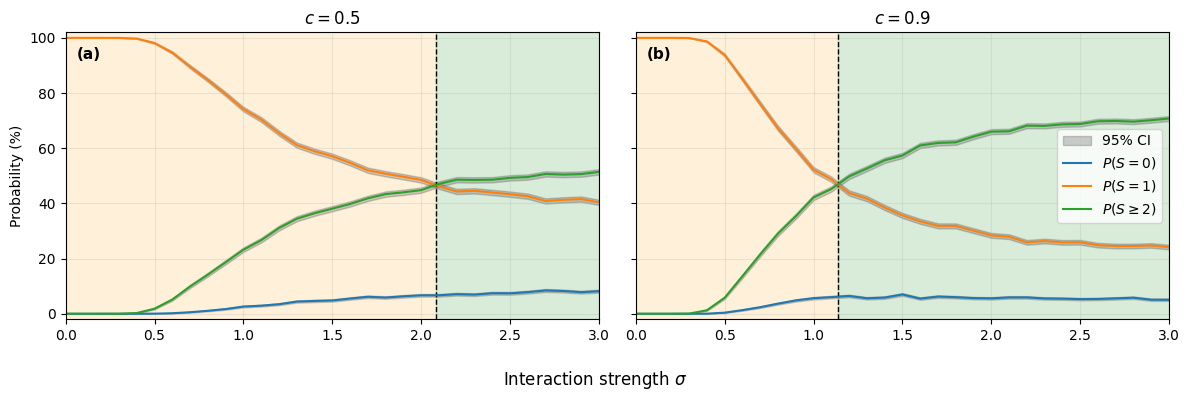}
    \caption{\textbf{Multiplicity regimes of stable states as a function of interaction strength $\sigma$ in the structured GLV model.}
    The curves show the probabilities of observing no stable state, one stable state, or multiple stable states. Same parameter sweep, confidence interval convention and dominance region convention as in Fig.~\ref{fig:glv_sat_mono_multi}, but computed for stable states.
    }
    \label{fig:glv_sta_mono_multi}
\end{figure}

The stability filter effect is further reflected in the multiplicity regimes shown in Fig.~\ref{fig:glv_sta_mono_multi}.

As in the saturated state case, increasing $\sigma$ drives a transition from a regime dominated by a single state to a regime where multiple states become more likely, and higher connectance shifts this transition toward smaller values of $\sigma$. However, this transition is less pronounced for stable states. In particular, $P(S\geq2)$ increases more slowly, while $P(S=0)$ becomes visible at intermediate to high interaction strengths. Stability therefore counterbalances the multiplicity generated at the saturated level: stronger and denser interactions may create several saturated states, but only a subset of them remains dynamically stable. 

The size component tells a complementary story. 
The regime in which stable states become more numerous is not necessarily a regime in which they become larger.

\begin{figure}[!htbp]
    \centering
    \includegraphics[width=0.7\linewidth]{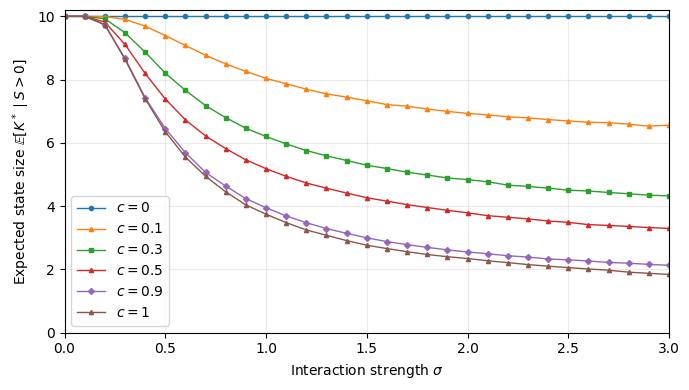}
    \caption{\textbf{Expected number of coexisting species conditionally on the existence of stable states as a function of interaction strength $\sigma$ in the structured GLV model for fixed connectance values.}
    The curves show $\mathbb{E}[K^\ast\mid S>0]$ for fixed connectance values. Same parameter sweep and sampling procedure as in Fig.~\ref{fig:glv_sat_exp}, but computed for stable states.
    }
    \label{fig:glv_sta_exp}
\end{figure}

The expected size of stable states, shown in Fig.~\ref{fig:glv_sta_exp}, provides complementary information. As for saturated states, the disconnected case $c=0$ remains at the maximal value $\mathbb{E}[K^\ast\mid S>0]=N$, and all curves start from this value at $\sigma=0$. Once interactions are introduced, the expected size decreases more strongly for stable states than for saturated states. Increasing connectance amplifies this decrease: denser interaction networks make large stable coexistence states harder to maintain, leading to lower typical number of coexisting species at high interaction strength.

Together with Fig.~\ref{fig:glv_sta_mono_multi}, this reveals the central pattern of the stable GLV landscape: increasing interaction strength makes stable states more numerous within systems but smaller in support size.

Overall, the GLV example shows that the proposed formalism turns the equilibrium landscape into a set of system-level observables. 
Saturation describes the supports allowed by non-invadability, whereas stability shows how this landscape is filtered by internal dynamics. 
Across the structured ensemble, interaction strength and connectance affect not only the existence of admissible states, but also their multiplicity and their typical support size.

\section{Replicator dynamics}
\subsection{Model definition and equilibrium types}
The replicator equation describes the evolution of species frequencies 
$z=(z_1,\ldots,z_N)$ in the probability simplex
\[
\Sigma_N=\left\{z\in[0,1]^N \mid \sum_{i=1}^N z_i=1\right\},
\]
and is governed by a fitness matrix 
\(\Lambda=(\lambda_i^j)_{1\leq i,j\leq N}\) 
\citep{HofbaeurSigmund1998, Nowak2006}:
\begin{equation}
    \dot z_i = z_i\big((\Lambda z)_i - z^\top \Lambda z\big), \forall i \in \{1, \ldots, N\}.
    \label{eq:replicator}
\end{equation}

\begin{definition}[Admissible steady state]\label{def:s_replicator}
     $z^\ast \in \Sigma_N$ is an \emph{admissible steady state} if it satisfies
    \begin{equation}
        z_i^\ast \bigl( (\Lambda z^\ast)_i - (z^\ast)^\top \Lambda z^\ast \bigr) = 0,
        \quad \forall i = 1,\ldots,N.
        \label{eq:s_replicator}
    \end{equation}
\end{definition}

The admissible steady states are given by setting \eqref{eq:replicator} to zero in $\Sigma_N$. Since $z^\ast \in \Sigma_N$, at least one $z^\ast_i$ must be positive. The $2^N-1$ non-empty supports correspond to all possible subcommunities of the $N$ species considered. Since each support can generate at most one admissible steady state, the total number of admissible steady states is at most $2^N-1$.

\begin{definition}[Saturated state]\label{def:sat_replicator}
    A steady state $z^\ast \in \Sigma_N$ is said to be \emph{saturated} if
    \begin{equation}
        (\Lambda z^\ast)_i - (z^\ast)^\top \Lambda z^\ast \le 0,
        \quad \forall i = 1,\ldots,N.
        \label{eq:sat_replicator}
    \end{equation}
\end{definition}

The interpretation of saturated states extends directly to the replicator dynamics. The saturation condition \eqref{eq:sat_replicator} is precisely the criterion for a Nash equilibrium in game theory \citep{Nash1950, HofbaeurSigmund1998}.

\begin{definition}[Asymptotically stable state]\label{def:sta_replicator}
    A steady state $z^\ast \in \Sigma_N$ is said to be \emph{asymptotically stable} if all eigenvalues of its Jacobian matrix restricted to the tangent space of $\Sigma_N$ have strictly negative real parts.
\end{definition}

As for the GLV model, the Jacobian admits a block structure separating internal dynamics on the support from invasion directions, implying that saturation is a necessary condition for stability.

Beyond mere dynamical stability, we consider a stronger notion of robustness which is evolutionary stability \citep{HofbaeurSigmund1998}. 

\begin{definition}[evolutionarily stable state (ESS)]\label{def:ESS}
    A steady state $z^\ast \in \Sigma_N$ is said to be an \emph{ESS} if
    \begin{enumerate}
        \item $z^\ast$ is a saturated state, i.e. satisfies the inequalities \eqref{eq:sat_replicator},
        \item and if equality holds for some $p \ne z^\ast$, then $p^\top \Lambda p < {z^\ast} ^\top \Lambda p$.
    \end{enumerate}
\end{definition}

Starting from Definition~\ref{def:ESS}, one can derive an explicit characterization of ESS that is particularly convenient for numerical implementation. Let $z^\ast \in \Sigma_N$ denote a steady state with support $\mathrm{supp}(z^\ast)$. Then $z^\ast$ is an ESS if and only if it satisfies the quadratic inequality
\begin{equation}
    (z^\ast - p)^\top \Lambda (z^\ast - p) < 0 \quad \forall\,  p \in \Sigma_N, \, p \neq  z^\ast \, \text{ s.t } \, \mathrm{supp}(p) \subseteq \mathrm{supp}(z^\ast).
    \label{eq:ESS_quadratic}
\end{equation}

This condition expresses the local non-invadability of $z^\ast$ against all nearby mutant frequency distributions supported on the same set of species \citep{BroomRychtar2013}.
Under these conditions, an ESS is globally stable relative to the subsimplex defined by its support \citep{HofbaeurSigmund1998}.

These notions form the hierarchy $\mathrm{ESS} \subset \text{asymptotically stable}  \subset \mathrm{saturated} \subset \mathrm{steady}$, and delimit the state set $E_p$ on which multiplicity and biodiversity are evaluated.

\subsection{Interaction structures and numerical procedure}

This section specifies how the parameters are generated and how the equilibrium landscapes are computed for each realization.
Similarly to the GLV model, we analyze two relevant regimes for the fitness matrix $\Lambda$.
\paragraph{Random. }In the first configuration, without loss of generality, the fitness matrix $\Lambda$ is defined by $\lambda_i^i = 0$, and off-diagonal entries $\lambda_i^j$ correspond to pairwise invasion fitness and are independently drawn sampled from a centered Gaussian distribution with unit variance. This structure represents a generic interaction landscape with no built-in biological asymmetry beyond randomness.
\paragraph{Structured ($\mu$). }In the second configuration, the fitness matrix $\Lambda(\mu)$ is derived from a multi-strain SIS co-colonization model \citep{MadecGjini2021}. In this setting, the invasion fitness of strain $i$ into a population dominated by strain $j$ is given by
\begin{equation}
    \lambda_i^j = \alpha_{j,i} - \alpha_{j,j} + \mu \bigl( \alpha_{j,i} - \alpha_{i,j} \bigr),
    \label{eq:epi_lambda}
\end{equation}
where $A = (\alpha_{i,j})_{1 \le i,j \le N}$ encodes pairwise co-colonization interactions. The coefficients $\alpha_{i,j}$ are drawn independently from a centered Gaussian distribution with unit variance.

\begin{table}[H]
\centering
\caption{\textbf{Summary of the fitness matrix structures considered in the replicator model.}}
\begin{tabular}{|l|l|}
\hline
Structure & Fitness matrix $\Lambda$ \\
\hline
Random 
& $\lambda_i^i = 0$, \quad 
$\lambda_i^j \sim \mathcal{N}(0,1)$ for $i \neq j$ \\[4pt]

\hline

Structured ($\mu$)
& $\lambda_i^j = \alpha_{j,i} - \alpha_{j,j}
+ \mu(\alpha_{j,i} - \alpha_{i,j})$, \\
& $\alpha_{i,j} \sim \mathcal{N}(0,1)$ \\

\hline
\end{tabular}

\label{tab:rep_structures}
\end{table}

For each realizations, all $2^N - 1$ non-empty supports are enumerated. For each support, the corresponding steady state is computed by solving \eqref{eq:s_replicator}, generically yielding a finite set of admissible steady states.

The resulting states are first tested for saturation \eqref{eq:sat_replicator}, then for stability, and finally for evolutionary stability using the quadratic criterion \eqref{eq:ESS_quadratic}. Classifying the states hierarchically reduces the computational cost; for instance, no need to compute the Jacobian for the non-saturated states.

\subsection{Results}
\label{sec:rep_results}

Now we apply the same formalism to the replicator model. 
The aim is not to introduce new quantities, but to test whether the same decomposition, i.e. existence, multiplicity, and state size composition, remains informative in a different dynamical setting. 
For $N=10$ and each interaction structure, we generate $M = 10\,000$ independent systems, enumerate all candidate supports ($2^N-1$), and classify the corresponding equilibria as saturated, asymptotically stable, or evolutionarily stable.

The analysis follows the same logic as in the GLV case.

\subsubsection{Saturated states}
\label{sec:rep_saturated_results}

We first consider saturated states. 
As in the GLV model, saturated states describe admissible states that are resistant to invasion, but now within the geometry of the simplex.

\begin{figure}[!htbp]
    \centering

    \begin{minipage}{\textwidth}
        \centering
        \text{(I) Random }\\[0.2em]
        \includegraphics[width=\textwidth]{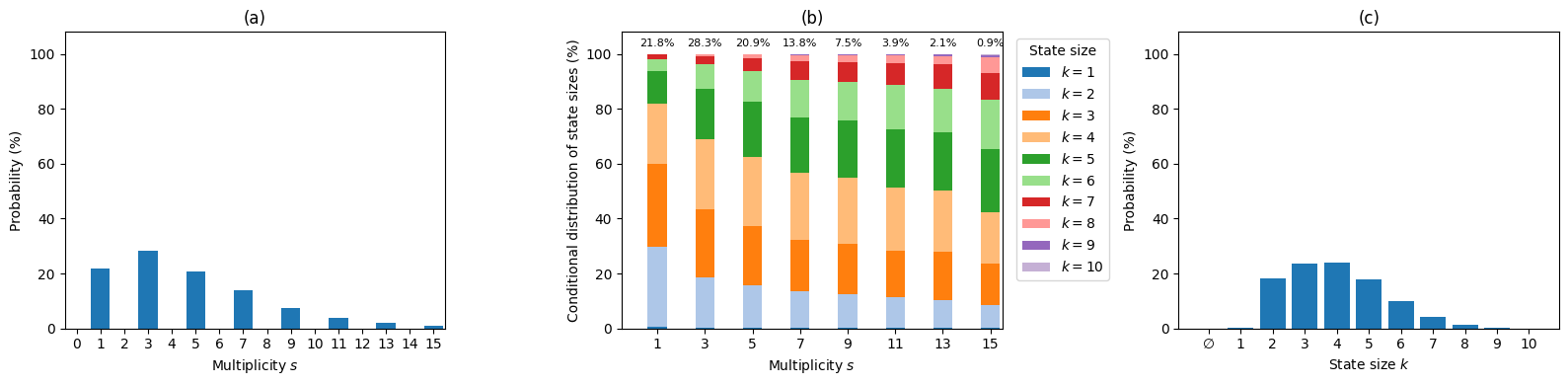}
    \end{minipage}

    \vspace{0.8em}

    \begin{minipage}{\textwidth}
        \centering
        \text{(II) $\mu=0.1$}\\[0.2em]
        \includegraphics[width=\textwidth]{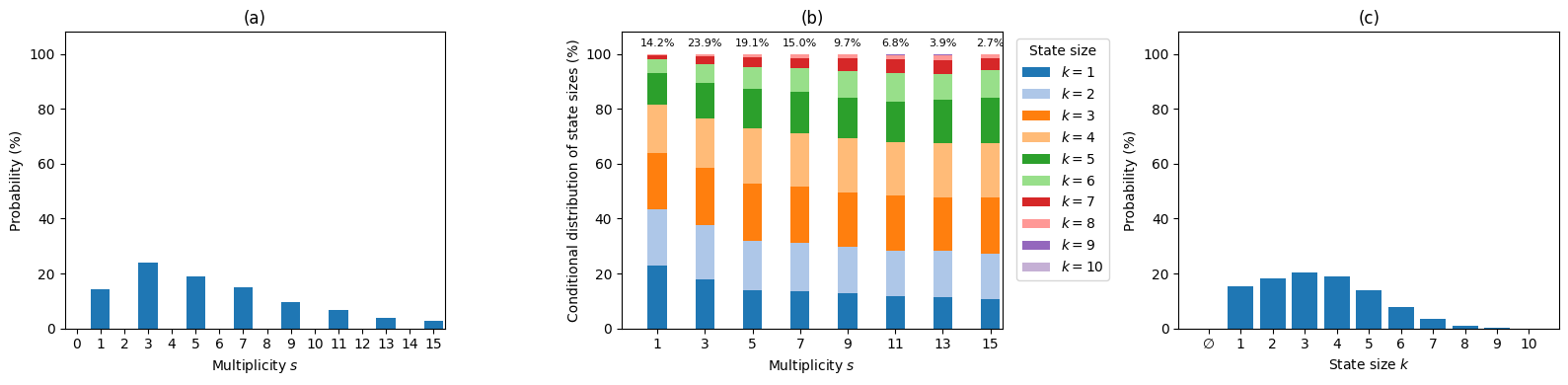}
    \end{minipage}

    \vspace{0.8em}

    \begin{minipage}{\textwidth}
        \centering
        \text{(III) $\mu=10$}\\[0.2em]
        \includegraphics[width=\textwidth]{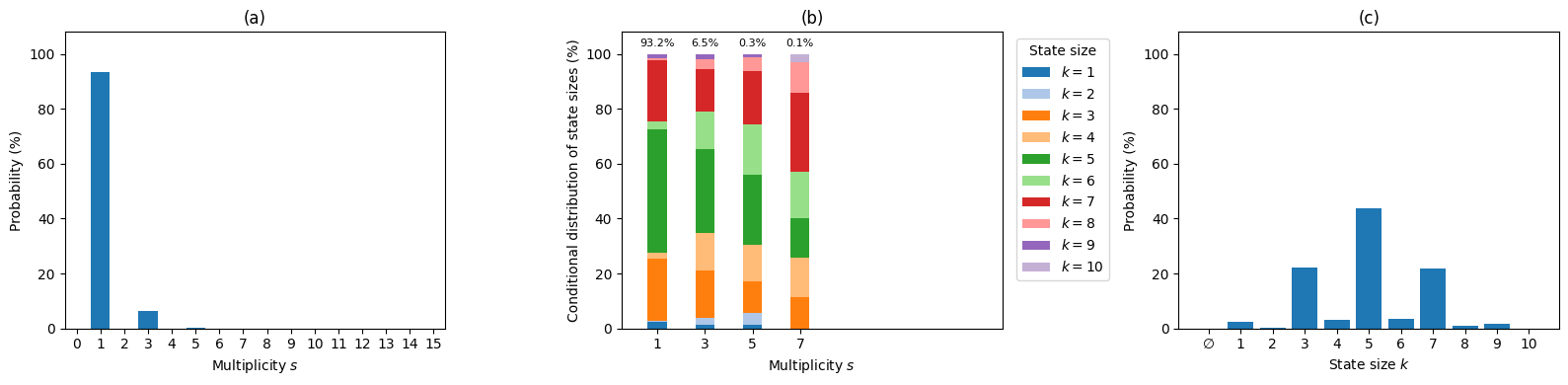}
    \end{minipage}

    \caption{
    \textbf{Output of the formalism for saturated states in the replicator model under different regimes.}
    The three rows correspond to (I) the random regime, (II) the structured regime with $\mu=0.1$, and (III) the structured regime with $\mu=10$, as defined in Table~\ref{tab:rep_structures}.
    Saturated states are identified by enumerating all possible supports.
    Within each row, Panel~(a) shows the distribution of the number $s$ of saturated states, Panel~(b) shows the conditional distribution of state sizes $k$ within each multiplicity class, and Panel~(c) shows the marginal distribution of state sizes aggregated over all systems and multiplicities.
    In Panel~(b), each stacked bar is normalized to $100\%$, while the percentages displayed above the bars indicate the fraction of systems belonging to the corresponding multiplicity class.
    For clarity, multiplicity classes are displayed up to $s_{\max}=15$.
    }
    \label{fig:rep_sat_formalism_regimes}
\end{figure}

As before, Fig.~\ref{fig:rep_sat_formalism_regimes} separates information that would otherwise be pooled together: the number of saturated states per system, the state size composition within each multiplicity class, and the aggregated size distribution.

The same properties are observed in the replicator model in the random and $\mu$-structured case, as in the GLV model: the existence, the occurrence of odd multiplicities (Panel~(a)). This similarity may be explained by the equivalence between the two models \citep{HofbaeurSigmund1998, Allesina2026}. From the perspective of game theory, these properties are consistent with those of Nash equilibria, namely their existence \citep{Nash1950} and the oddness of their number \citep{Harsanyi1973}.

In the random regime (I), the multiplicity distribution is mainly concentrated essentially between $s=1$ and $s=7$. 
Within each multiplicity class, the state size composition remains heterogeneous, but states generally involve at least two coexisting species; states with a single species are almost absent, indicating that they are not resistant to invasion. This is confirmed by the state size distribution in Panel~(c), which is unimodal and concentrated on intermediate support sizes. The mass at $k=1$ is negligible, the distribution effectively starts at $k=2$, and large support sizes become rare.

The structured case allows us to compare two contrasting interaction regimes: (II) $\mu=0.1$, where double infections dominate, and (III) $\mu=10$, which corresponds to the prevalence of single infections \citep{MadecGjini2021}. 

Comparing both cases shows that increasing $\mu$ reduces the multiplicity of saturated state multiplicity. Because stable states and ESS are obtained by imposing additional constraints on saturated states, this suggests that their multiplicities should follow the same trend.

At $\mu=0.1$, systems display a broader range of multiplicities, with a state size distribution spread over small and moderate values. Whereas at $\mu=10$, most systems have only one or a few saturated states, and the distribution of state sizes is concentrated on odd values. This is due to the asymptotic antisymmetric structure of the fitness matrix as $\mu \rightarrow \infty$ \citep{chawanya2002, MadecGjini2021}.

\begin{figure}[!htbp]
    \centering
    \includegraphics[width=0.7\linewidth]{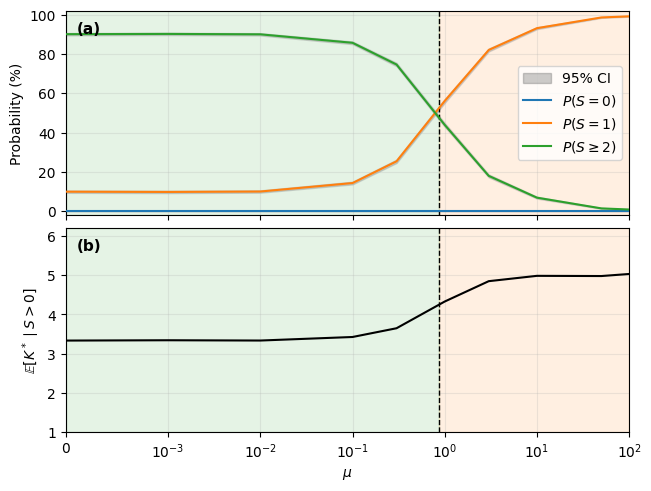}
    \caption{
    \textbf{Multiplicity regimes and conditional expected size of saturated states as functions of the ratio of single to co-colonization $\mu$ in the structured replicator model.}
    For each value of $\mu$ varied from $0$ to $100$, with denser sampling near small and intermediate values, we generate independent systems using the same sampling procedure as in Table~\ref{tab:rep_structures}, and classify each system according to the number $S$ of saturated states it admits.
    (a) Probabilities of observing no saturated state, $P(S=0)$, exactly one saturated state, $P(S=1)$, or multiple saturated states, $P(S\geq2)$. The grey bands represent approximate $95\%$ confidence intervals.
    (b) Conditional expected state size $\mathbb{E}[K^\ast\mid S>0]$.
    In both Panels, the shaded regions indicate the dominant multiplicity regime, and the vertical dashed line marks the intersection between $P(S=1)$ and $P(S\geq2)$.
    }
    \label{fig:rep_sat_mono_multi_exp}
\end{figure}

To determine whether these differences correspond to isolated regimes or to a systematic transition, we next vary $\mu$ while retaining the same sampling procedure (Fig.~\ref{fig:rep_sat_mono_multi_exp}).

At small $\mu$, multiple saturated states are more likely to occur than a unique saturated state, consistent with the broad multiplicity distribution observed for $\mu=0.1$. As $\mu$ increases, $P(S\geq2)$ decreases while $P(S=1)$ increases, and the two probabilities become equal at the transition marked by the vertical dashed line. Beyond this point, the uniqueness of saturated states rapidly becomes dominant, approaching near certainty for large $\mu$.

The conditional expected state size provides a complementary description of the same reorganization. For low values of $\mu$ in the multi-state regime, saturated states have a smaller expected size, whereas the progressive dominance of a unique saturated state is accompanied by a small increase in $\mathbb{E}[K^\ast\mid S>0]$. The most pronounced variation in expected size occurs around the same range of $\mu$ over which the multiplicity probabilities cross. Hence, the transition is not limited to a reduction in the number of saturated states, it is also associated with a restructuring of their composition toward slightly larger states.

\subsubsection{Stable states}
\label{sec:replicator_stable_results}

\begin{figure}[!htbp]
    \centering

    \begin{minipage}{\textwidth}
        \centering
        \text{(I) Random }\\[0.2em]
        \includegraphics[width=\textwidth]{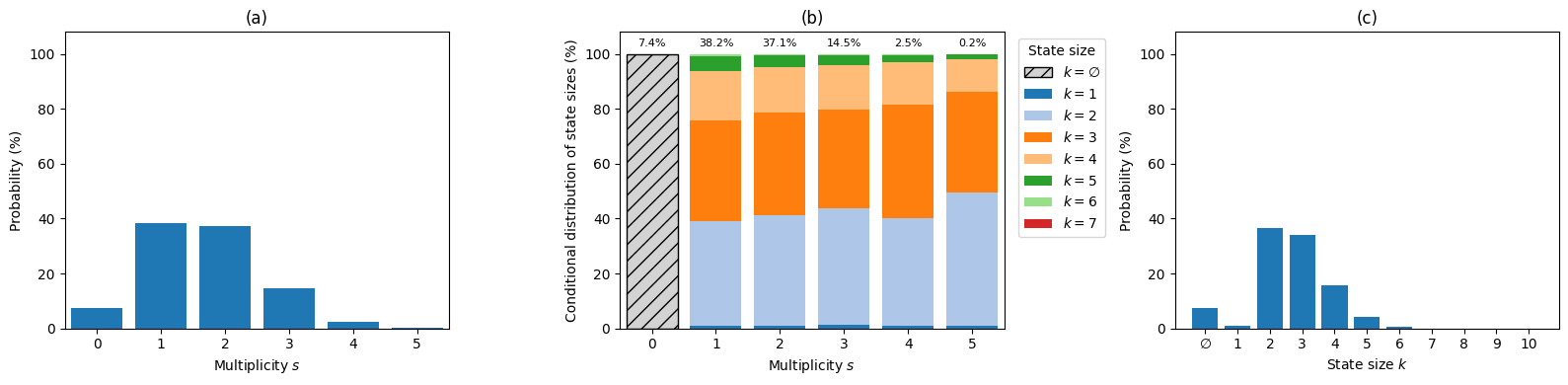}
    \end{minipage}

    \vspace{0.8em}

    \begin{minipage}{\textwidth}
        \centering
        \text{(II) $\mu=0.1$}\\[0.2em]
        \includegraphics[width=\textwidth]{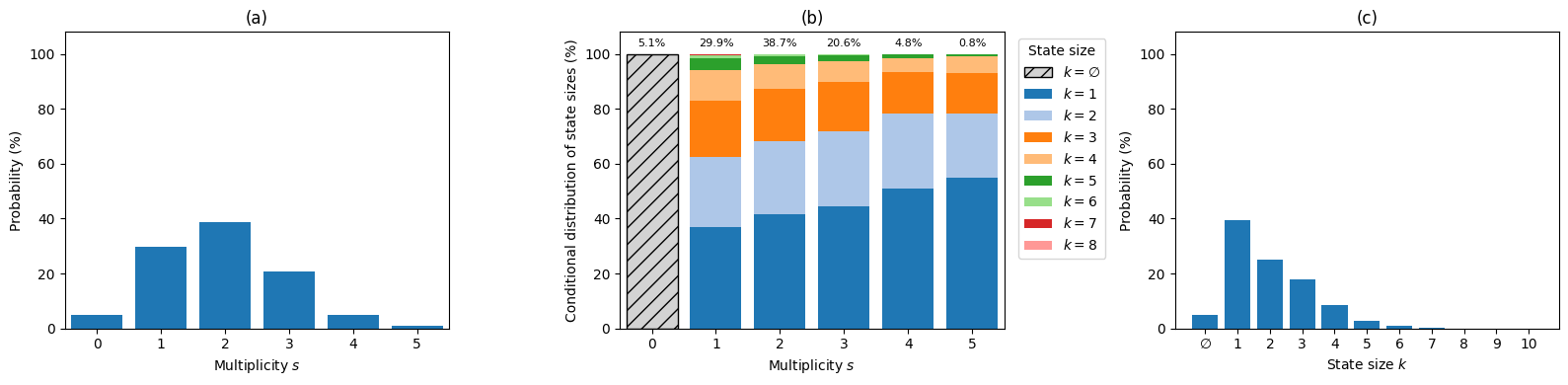}
    \end{minipage}

    \vspace{0.8em}

    \begin{minipage}{\textwidth}
        \centering
        \textbf{(III) $\mu=10$}\\[0.2em]
        \includegraphics[width=\textwidth]{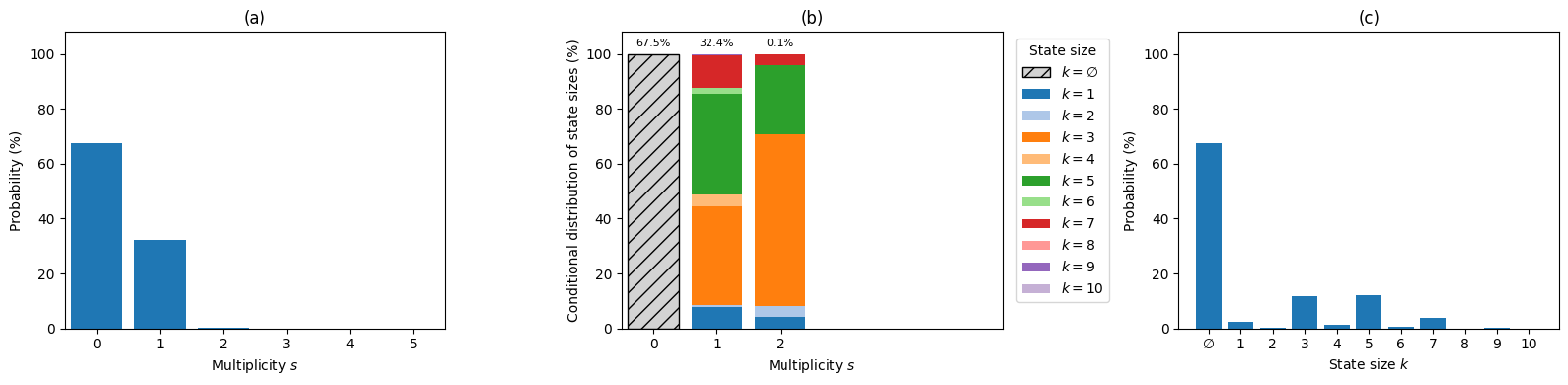}
    \end{minipage}

    \caption{
    \textbf{Output of the formalism for stable states in the replicator model under different regimes.}
    The three rows correspond to (I) the random regime, (II) the structured regime with $\mu=0.1$, and (III) the structured regime with $\mu=10$, as defined in Table~\ref{tab:rep_structures}.
    Within each row, Panel~(a) shows the distribution of the number $s$ of stable states, Panel~(b) shows the conditional distribution of state sizes $k$ within each multiplicity class, and Panel~(c) shows the marginal distribution of state sizes.
    The hatched $k=\varnothing$ category represents systems with no stable state.
    The sampling procedure and graphical conventions are the same as in Fig.~\ref{fig:rep_sat_formalism_regimes}, but the analysis is restricted to stable states.
    }
    \label{fig:rep_stable_formalism_regimes}
\end{figure}

We then impose local stability on the saturated states. 
This additional constraint substantially reorganizes the distributions shown in Fig.~\ref{fig:rep_stable_formalism_regimes}. 

In contrast to saturated states, the existence of a stable state is no longer guaranteed, and the oddness of the multiplicity property observed at the saturated level is
lost. Stable states may therefore be absent, and both odd and even multiplicities occur. Stability consequently acts as a nontrivial filter of the saturated state set rather than simply reducing its overall multiplicity. In the random regime, the multistability distribution is concentrated mainly around one and two states, with these two classes occurring with comparable probabilities and compressed toward smaller values once stability is imposed. The corresponding state size distribution is also shifted toward smaller supports, indicating that many of the larger saturated configurations do not persist as stable states. 
For $\mu=0.1$, most systems still admit stable states, and the multiplicity is mainly concentrated around small values of $s$. 
For $\mu=10$, absence becomes the dominant outcome: most systems have no stable state, indicating more complex dynamics, while the few systems that admit stable states are mostly monostable. 
This marked increase in the empty class is precisely the type of situation motivating the hurdle and zero-inflated structure of the formalism, where the probability of existence must be separated from the conditional distribution of state-sizes when they exist. 
In the structured replicator model, $\mu$ represents the ratio of single to co-colonization, a parameter known to strongly affect coexistence complexity in multi-strain systems \citep{MadecGjini2021}.

\begin{figure}[!htbp]
    \centering
    \includegraphics[width=0.7\linewidth]{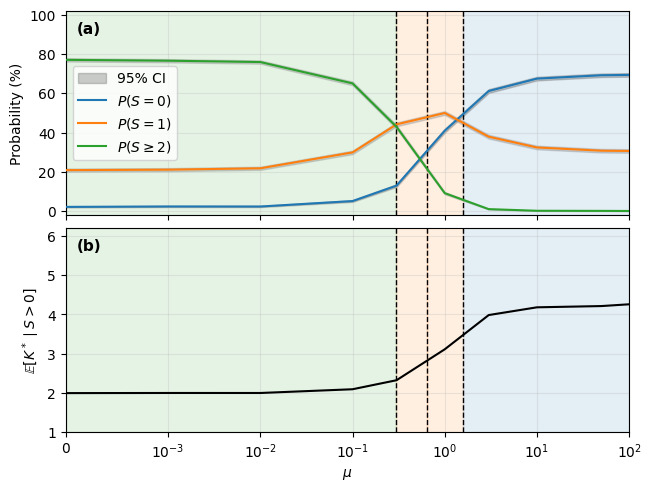}
    \caption{
    \textbf{Multiplicity regimes and conditional expected size of stable states as functions of the ratio of single to co-colonization $\mu$ in the structured replicator model.}
    (a) Probabilities of observing no stable state, $P(S=0)$, exactly one stable state, $P(S=1)$, or multiple stable states, $P(S\geq2)$.
    (b) Conditional expected stable-state size $\mathbb{E}[K^\ast\mid S>0]$.
    The parameter sweep, sampling procedure, confidence interval convention, and graphical conventions are the same as in Fig.~\ref{fig:rep_sat_mono_multi_exp}, but the analysis is restricted to stable states.
    Shaded regions indicate the dominant multiplicity regime, and the vertical dashed lines mark pairwise intersections that these regimes.
    }
    \label{fig:rep_stable_mono_multi_exp}
\end{figure}

Fig.~\ref{fig:rep_stable_mono_multi_exp} shows that $\mu$ reorganizes the stable state outcomes in two complementary ways. 
For small $\mu$, most systems admit multiple stable states, but these states are typically small, involving only about two coexisting species. 
As $\mu$ increases, multiplicity collapses. At the same time, the stable states that remain become larger on average. 
Thus, increasing $\mu$ shifts the system from many small alternative stable states to fewer, but larger stable coexistence states. 
The dominance of $\mathbb{P}(S=0)$ for large values of $\mu$ suggests a shift toward other types of dynamics including limit cycles and chaos. 
 
This behavior is consistent with recent studies emphasizing that multistability in complex communities should be characterized not only by the number of stable states, but also by their richness and composition \citep{Guimera2024,Aguade2024}.

\subsubsection{Evolutionarily stable states}
\label{sec:replicator_ess_results}

The replicator model also allows us to consider a third class of states. 
ESSs provide a stronger selection criterion than feasibility, saturation, or asymptotical stability alone. 

Therefore, we apply the same system-centered formalism to this more restrictive class of states.

\begin{figure}[!htbp]
    \centering

    \begin{minipage}{\textwidth}
        \centering
        \textbf{(I) Random}\\[0.2em]
        \includegraphics[width=\textwidth]{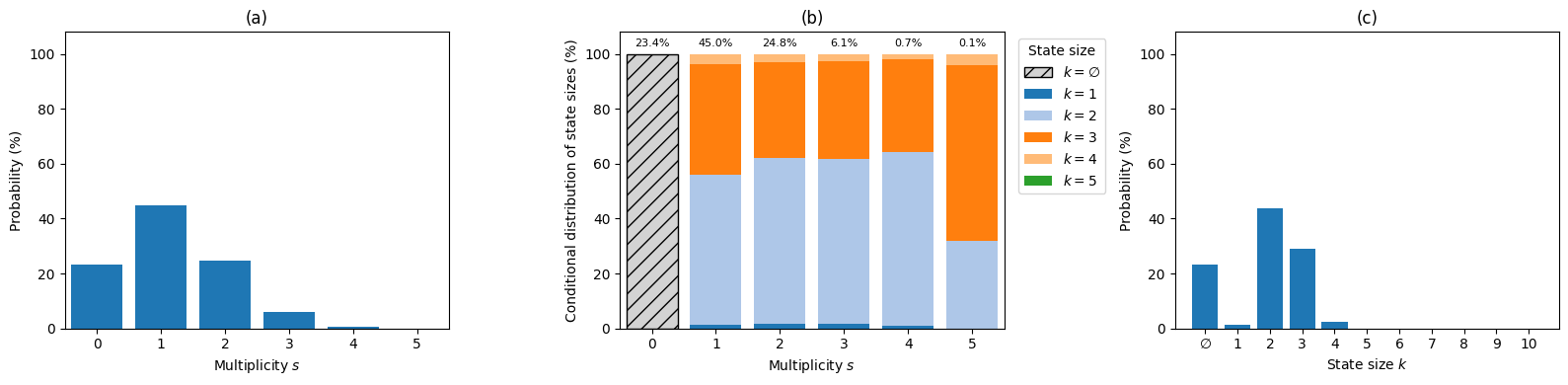}
    \end{minipage}

    \vspace{0.8em}

    \begin{minipage}{\textwidth}
        \centering
        \textbf{(II) $\mu=0.1$}\\[0.2em]
        \includegraphics[width=\textwidth]{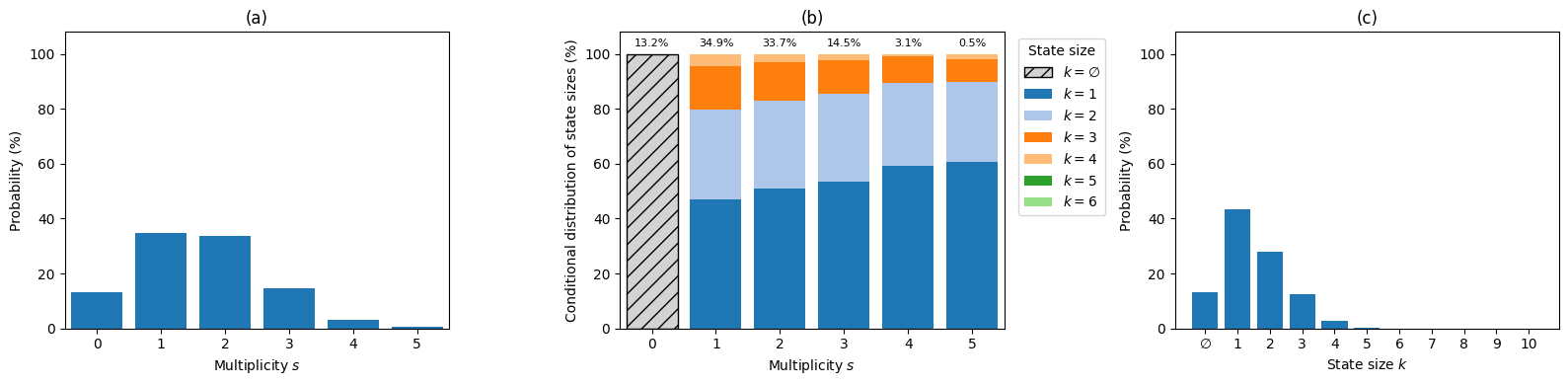}
    \end{minipage}

    \vspace{0.8em}

    \begin{minipage}{\textwidth}
        \centering
        \textbf{(III) $\mu=10$}\\[0.2em]
        \includegraphics[width=\textwidth]{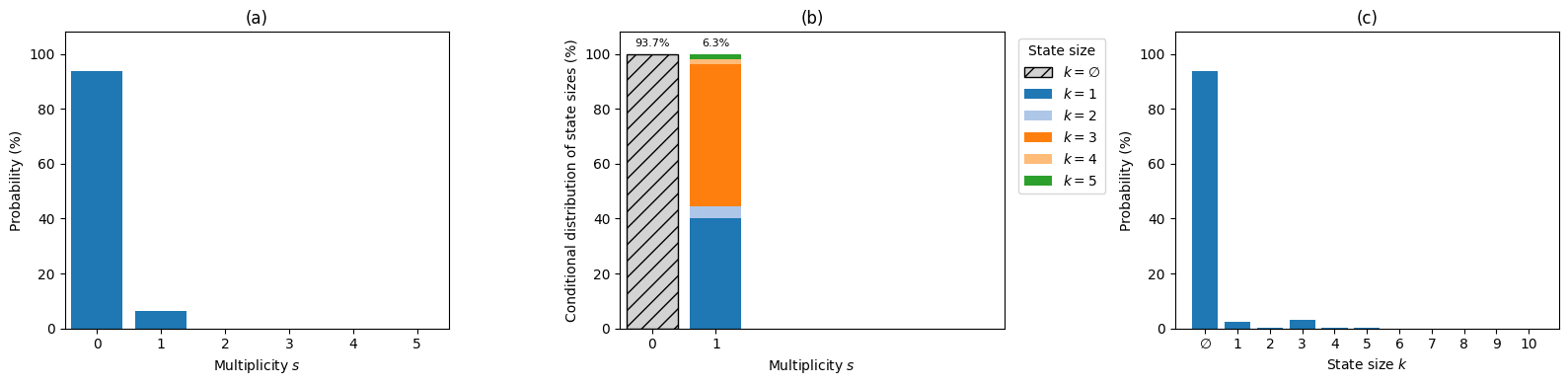}
    \end{minipage}

    \caption{
    \textbf{Output of the formalism for evolutionarily stable states in the replicator model under different regimes.}
    The three rows correspond to (I) the random regime, (II) the structured regime with $\mu=0.1$, and (III) the structured regime with $\mu=10$, as defined in Table~\ref{tab:rep_structures}.
    The sampling procedure and graphical conventions are the same as in Fig.~\ref{fig:rep_stable_formalism_regimes}.
    }
    \label{fig:rep_ess_formalism_regimes}
\end{figure}

Fig.~\ref{fig:rep_ess_formalism_regimes} reveals the strong filtering effect of the ESS criterion. 
The random ensemble and the $\mu=0.1$ regime still admit ESS in a non-negligible fraction of systems, mainly with low multiplicity and small state sizes. 
At $\mu=10$, however, the empty class dominates, and the few remaining ESS are almost always unique. 
Thus, increasing $\mu$ primarily reduces the existence of ESS states, highlighting the importance of separating absence from the conditional distribution of states when they exist.

\begin{figure}[!htbp]
    \centering
    \includegraphics[width=0.7\linewidth]{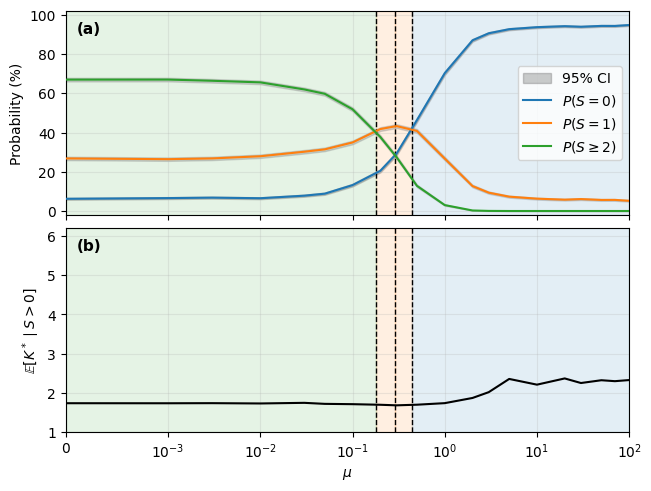}
    \caption{
    \textbf{Multiplicity regimes and conditional expected size of evolutionarily stable states as functions of the ratio of single to co-colonization $\mu$ in the structured replicator model.}
    (a) Probabilities of observing no ESS, $P(S=0)$, exactly one ESS, $P(S=1)$, or multiple ESSs, $P(S\geq2)$.
    (b) Conditional expected ESS size $\mathbb{E}[K^\ast\mid S>0]$.
    The parameter sweep, sampling procedure, confidence interval convention, and graphical conventions are the same as in Fig.~\ref{fig:rep_stable_mono_multi_exp}.
    Shaded regions indicate the dominant multiplicity regime, while vertical dashed lines mark pairwise intersections between the probability curves.
    }
    \label{fig:rep_ess_mono_multi_exp}
\end{figure}

Fig.~\ref{fig:rep_ess_mono_multi_exp} reveals a strong separation between existence and state-size distribution. 
At low $\mu$, systems are predominantly characterized by multiple ESSs. This multistable regime is followed by a narrow intermediate range in which a unique ESS is the most likely outcome. At larger $\mu$, the probability of finding no ESS increases rapidly and becomes overwhelmingly dominant.The comparison with Fig.~\ref{fig:rep_stable_mono_multi_exp} highlights the additional selectivity imposed by evolutionary stability. Both transitions occur at smaller values of $\mu$ for ESSs, and the intermediate regime dominated by a unique state is considerably narrower. In particular, the absence of ESSs becomes dominant while stable states are still frequently observed. Thus, increasing $\mu$ does not merely reduce multiplicity; it rapidly eliminates evolutionary robustness, even when dynamically stable states remain available. The conditional expected size provides a complementary distinction between the two state classes. While the expected size of the stable states increases substantially with $\mu$, the expected ESS size remains close to the small supports and rises only moderately after the transition to absence. Consequently, the rare ESSs that persist at large $\mu$ remain much smaller than the corresponding stable states. Together, these results show that the ESS condition filters both the existence and the composition of states.

\begin{figure}[!htbp]
    \centering
    \includegraphics[width=0.7\linewidth]{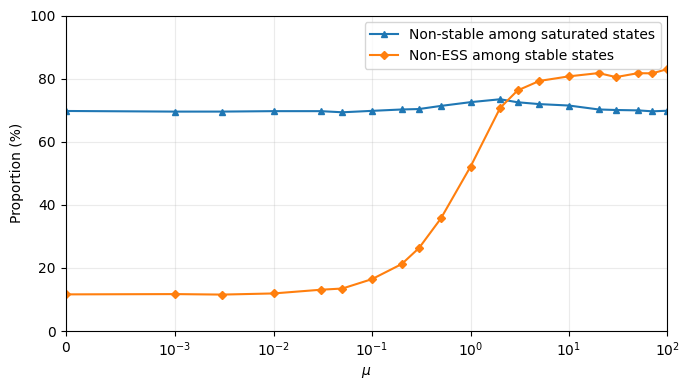}
    \caption{
    \textbf{Filtering between saturated, stable and evolutionarily stable states in the structured replicator model.}
    The blue curve shows the proportion of saturated states that are not stable, and the orange curve shows the proportion of stable states that are not ESS, as functions of $\mu$.
    Proportions are computed from expected multiplicities.
    }
    \label{fig:rep_filtering_proportions}
\end{figure}

To make the successive filtering effects more explicit, Fig.~\ref{fig:rep_filtering_proportions} compares the loss of states between the three levels of states. For each value of $\mu$, we compute the proportion of saturated states that fail to be stable and the proportion of stable states that fail to be ESS. These quantities are obtained from the expected multiplicities as
\[
\frac{\mathbb{E}[S^{\mathrm{saturated}}]-\mathbb{E}[S^{\mathrm{stable}}]}
{\mathbb{E}[S^{\mathrm{saturated}}]}
\qquad\text{and}\qquad
\frac{\mathbb{E}[S^{\mathrm{stable}}]-\mathbb{E}[S^{\mathrm{ESS}}]}
{\mathbb{E}[S^{\mathrm{stable}}]}.
\]

The stability constraint removes a large fraction of the saturated states across the entire parameter range. This proportion remains close to $70\%$ and varies only weakly with $\mu$, indicating that stability imposes a strong constraint independently of the ratio $\mu$. On the contrary, the evolutionary stability filter, from stable states to ESSs, is strongly dependent of $\mu$. For small $\mu$, only a limited fraction of stable states fails to be ESS, whereas this fraction increases sharply around intermediate values of $\mu$ and exceeds $80\%$ for large $\mu$.

\section{Discussion}
\label{sec:discussion}

The central idea of this work is to take the system, rather than the individual state, as the primary statistical unit. 
For each sampled parameter set, the model defines a system. The states associated with that system are then enumerated, classified, and summarized through their existence, multiplicity, and size composition. 
This system-centered perspective complements the usual state-centered approach, in which all states are pooled across systems before computing summary statistics.

\subsection{System-centred versus state-centred descriptions}

The formalism separates three levels of information that are often mixed together: whether a system admits at least one state of a certain type, how many such states it admits, and how rich these states are. 
This decomposition is close in spirit to zero-inflated and hurdle-type models, where the probability of observing a nonzero outcome is separated from the conditional distribution of positive outcomes \citep{MULLAHY1986341, Lambert1992}. 
In the present setting, the first question is whether a system admits at least one state of a certain type. 
Conditionally on this event, one can then study the multiplicity of states and the distribution of their support sizes.

This distinction is particularly important in multistable systems. 
Multistability is not only the qualitative statement that several stable outcomes may exist. 
It also raises quantitative questions: How often do multiple states occur across parameter sets, how many states does a typical system admit, and how different are these states in size or composition? 
The multiplicity variable $S_P$ turns these questions into measurable observables. 

The state-centered approach remains natural and useful when the question of interest concerns a state drawn uniformly from the pooled collection of all admissible states. 
This perspective is consistent with the approaches of theoretical ecology that study the statistical properties of feasible or assembled equilibrium states in large ecological systems \citep{Bunin2017, Barbier2018, grilli2017feasibility}. 

Still, it answers a different sampling question. 
The state-centered distribution describes a typical state in the pooled ensemble. 
The system-centered distribution describes the typical state structure associated with a randomly sampled system.

The covariance term \eqref{eq:cov} makes this distinction explicit. 
The difference between the state-centered and system-centered mean state sizes is a size-biasing effect: weighting systems by their number of admissible states shifts the mean by a covariance term between this weight and the within-system mean state size, a standard identity behind weighted distributions \citep{patil}.

When each system admits exactly one state, $S_p=1$ for all $p$, and therefore
\(\operatorname{Cov}(S_p, F_p^k)=0\), the state-centered and system-centered approaches coincide, because pooling over states is equivalent to averaging over systems.
When this covariance is positive or negative, it quantifies the coupling between state multiplicity and state size.

The added value of the system-centered formalism is that it preserves information about how states are grouped by system. 
This information is lost when all states are pooled before analysis. 
In a pooled description, a system with many states automatically contributes more than a system with one state, and a system with no state contributes nothing. 
As a consequence, pooled observables cannot distinguish between an outcome based on many systems admitting one state and one dominated by fewer systems with many states.

\subsection{Multistability and coexistence in complex systems}

The examples studied here illustrate how the formalism quantifies multistability in concrete dynamical models. 
In the GLV model, saturated states describe admissible equilibria that cannot be invaded by absent species, whereas stable states are obtained by adding the internal Jacobian stability condition. 
The structured GLV ensemble shows that interaction strength and connectance do not only affect whether admissible or stable states exist. They also affect how many such states a system admits and how many species coexist within each state.

One of the main patterns observed in the GLV example is that stronger interactions can shift the landscape toward systems with more states, but with fewer coexisting species in each state. 
A regime with many alternative states is not necessarily a regime with high coexistence within each state. 
Conversely, a regime with large coexistence states may be dominated by uniqueness rather than multiplicity.
This perspective complements the classical complexity-stability theory, which focuses on the stability of equilibria or community matrices as a function of system size, connectance, and interaction strength \citep{May1972,AllesinaTang2012,Stone1988}. 

It also connects with more recent approaches to random ecological communities and high-dimensional GLV systems, where feasibility, coexistence, and multiple equilibria are studied as emergent properties of random interaction ensembles \citep{Bunin2017,Barbier2018}. 
Recent work on alternative stable states and coexisting subsets in complex ecosystems further emphasizes that the relevant object is often not a single equilibrium, but a collection of possible stable communities \citep{Guimera2024,Aguade2024}. 

The replicator examples show that the same construction is not specific to GLV systems. 
Although the state space, equilibrium conditions differ, the same system-centered quantities can be defined once each parameter set is associated with a finite collection of states. 
Saturated, stable, and evolutionarily stable states define different levels of selection within the same formal structure. 

\subsection{Generality of the framework}

The two examples considered in this paper should be interpreted as illustrations of the formalism rather than as its domain of validity. 
The framework is not specific to GLV systems, replicator dynamics, or to ecological interpretations of states. 
It applies whenever a parameter set defines a system and this system can be associated with a finite collection of states that can be classified according to a prescribed criterion.
The type of a state may vary from one model to another, but once this finite collection is defined, the same observables can be computed: probability of existence, multiplicity, and state-size distribution. 
The competitive GLV and replicator models are therefore examples of use, not limitations of the method.

For models with infinitely many admissible configurations, additional selection, discretization, or classification rules would be required before applying the present framework. This reflects a standard issue in dynamical systems theory: when the attractor is a continuum, counting states is no longer well defined without specifying an equivalence relation, a section of the state space, or another rule selecting representative states \citep{Strogatz2018}.
The finite state set assumption should be viewed as a structural condition of the current formalism. 

\subsection{Beyond state size: diversity across states}

We measure state size by the number of species present in a state \citep{Magurran2004}. 
This is the most direct notion of diversity and allows GLV and replicator systems to be compared within the same formalism. 
However, once a system may admit several states, instead of reducing biodiversity to the richness of each state considered separately, we can consider other diversity measures.

For example, if system $p$ admits states with supports
\[
I_1,\dots,I_{S_p},
\]
one may define the union diversity
\[
U_p =
\left|
\bigcup_{\alpha=1}^{S_p} I_\alpha
\right|,
\]
which counts the number of species that appear in at least one admissible state of the system. 
This quantity can be large even if each individual state is small, provided that different states involve different subsets of species. 
Thus, the framework can be extended from richness within states to diversity across states.

This distinction is important in multistable regimes. 
Biodiversity is then distributed not only within equilibria, but also across alternative equilibria. 
A system may have low coexistence within each state but high diversity across its possible states. 
The formalism therefore opens the possibility of studying biodiversity as a property of the whole state set, rather than only as a property of a single equilibrium.

\subsection{Limitations and outlook}

Our framework focuses on stationary states that can be enumerated from candidate supports, and does not take into consideration transient community compositions or non-stationary attractors such as cycles or chaos.
The implementation relies on enumerating all nonempty supports for replicator dynamics and all supports including extinction for GLV, then testing steadiness, saturation, stability, and, for replicator, ESS in a hierarchical way. This is exact for moderate $N$ but becomes prohibitive at high dimension. Here, we used \(N=10\), for which all \(2^{10}=1024\) supports can be systematically explored for each sampled system.

Several strategies can improve this approach. The problem is  parallel: different supports and different parameter realizations can be explored independently. Optimization-based approaches constrained by non-negativity and necessary stability criteria \citep{LISCHKE201724}, can drastically improve the computational time. For replicator dynamics, adding constraints according to the properties of ESS can significantly reduce the set of stable states.

Beyond computational limitations, the formalism defines a natural bridge between probabilistic assumptions on systems and probability distributions on their state structures. If parameters are sampled from a given distribution, the formalism turns this randomness into distributions of occurrence, multiplicity, and coexistence size. In the future, we could therefore seek analytical results connecting these system-centered summaries to properties of random matrices, using tools from random matrix theory and probability theory. Such results would turn the Monte Carlo summaries developed in this paper into theoretical results.
 
\section{Conclusion}

By jointly quantifying multiplicity and biodiversity, the present formalism provides a new way to characterize ecological and evolutionary systems. Rather than treating equilibria as isolated objects, the framework keeps track of the set of states generated by each system. It then asks whether a given type of equilibrium exists, how many alternative coexistence configurations a system can support and how coexistence size varies across these configurations.

This shift in perspective has implications for ecology, epidemiology, and evolutionary theory. Systems with similar average biodiversity may differ strongly in multiplicity, and thus in their sensitivity to perturbations, invasions, or environmental change.

The formalism is not limited to the two models considered, GLV and replicator. It can be applied to type of models for which the set of admissible states is finite.

It also opens the door to studying how control parameters such as $\mu$, $\sigma$ or $c$ reshape the entire outputs.

Possible extensions include alternative interaction structures, or larger system sizes. The approach could also be extended to structured interaction networks, temporal variability, or combined with basin of attraction analyses to bridge structural and dynamical notions of multistability.

\bibliographystyle{plainnat}
\bibliography{references}

\begin{thebibliography}{40}
\providecommand{\natexlab}[1]{#1}
\providecommand{\url}[1]{\texttt{#1}}
\expandafter\ifx\csname urlstyle\endcsname\relax
  \providecommand{\doi}[1]{doi: #1}\else
  \providecommand{\doi}{doi: \begingroup \urlstyle{rm}\Url}\fi

\bibitem[Aguadé-Gorgorió et~al.(2024)Aguadé-Gorgorió, Arnoldi, Barbier, and
  Kéfi]{Guimera2024}
G.~Aguadé-Gorgorió, J.-F. Arnoldi, M.~Barbier, and S.~Kéfi.
\newblock A taxonomy of multiple stable states in complex ecological
  communities.
\newblock \emph{Ecology Letters}, 27:\penalty0 e14413, 2024.
\newblock \doi{10.1111/ele.14413}.

\bibitem[Allesina(2026)]{Allesina2026}
Stefano Allesina.
\newblock Global stability of ecological and evolutionary dynamics via
  equivalence.
\newblock \emph{Proceedings of the National Academy of Sciences}, 123\penalty0
  (13):\penalty0 e2534915123, 2026.
\newblock \doi{10.1073/pnas.2534915123}.
\newblock URL \url{https://www.pnas.org/doi/abs/10.1073/pnas.2534915123}.

\bibitem[Allesina and Tang(2012)]{AllesinaTang2012}
Stefano Allesina and Si~Tang.
\newblock Stability criteria for complex ecosystems.
\newblock \emph{Nature}, 483:\penalty0 205--208, 2012.
\newblock \doi{10.1038/nature10832}.

\bibitem[Barbier et~al.(2018)Barbier, Arnoldi, Bunin, and Loreau]{Barbier2018}
Matthieu Barbier, Jean-François Arnoldi, Guy Bunin, and Michel Loreau.
\newblock Generic assembly patterns in complex ecological communities.
\newblock \emph{Proceedings of the National Academy of Sciences}, 115\penalty0
  (9):\penalty0 2156--2161, 2018.
\newblock \doi{10.1073/pnas.1710352115}.

\bibitem[Beisner et~al.(2003)Beisner, Haydon, and Cuddington]{Beisner2003}
Beatrix~E. Beisner, Daniel~T. Haydon, and Kim Cuddington.
\newblock Alternative stable states in ecology.
\newblock \emph{Frontiers in Ecology and the Environment}, 1\penalty0
  (7):\penalty0 376--382, 2003.

\bibitem[Broom and Rycht{\'a}r(2013)]{BroomRychtar2013}
Mark Broom and Jan Rycht{\'a}r.
\newblock \emph{Game-Theoretical Models in Biology}.
\newblock CRC Press, 2013.
\newblock ISBN 9781003024682.

\bibitem[Bunin(2017)]{Bunin2017}
Guy Bunin.
\newblock Ecological communities with lotka-volterra dynamics.
\newblock \emph{Phys. Rev. E}, 95:\penalty0 042414, 2017.
\newblock \doi{10.1103/PhysRevE.95.042414}.

\bibitem[Chawanya and Tokita(2002)]{chawanya2002}
Tsuyoshi Chawanya and Kei Tokita.
\newblock Large-dimensional replicator equations with antisymmetric random
  interactions.
\newblock \emph{Journal of the Physical Society of Japan}, 71\penalty0
  (2):\penalty0 429--431, 2002.
\newblock \doi{10.1143/JPSJ.71.429}.

\bibitem[Clenet(2022)]{clenet}
Maxime Clenet.
\newblock \emph{{Large Lotka-Volterra model : when random matrix theory meets
  theoretical ecology}}.
\newblock Theses, {Universit{\'e} Gustave Eiffel}, December 2022.
\newblock URL \url{https://theses.hal.science/tel-04048703}.

\bibitem[Fujita et~al.(2025)Fujita, Yoshida, Suzuki, and Toju]{article}
Hiroaki Fujita, Shigenobu Yoshida, Kenta Suzuki, and Hirokazu Toju.
\newblock Alternative stable states of microbiome structure and soil ecosystem
  functions.
\newblock \emph{Environmental Microbiome}, 20, 03 2025.
\newblock \doi{10.1186/s40793-025-00688-4}.

\bibitem[Goh and Jennings(1977)]{GOH197763}
B.S. Goh and L.S. Jennings.
\newblock Feasibility and stability in randomly assembled lotka-volterra
  models.
\newblock \emph{Ecological Modelling}, 3\penalty0 (1):\penalty0 63--71, 1977.
\newblock ISSN 0304-3800.
\newblock \doi{https://doi.org/10.1016/0304-3800(77)90024-2}.
\newblock URL
  \url{https://www.sciencedirect.com/science/article/pii/0304380077900242}.

\bibitem[Gore et~al.(2025)Gore, Hu, He, Barbier, Song, and Bunin]{unknown}
Jeff Gore, Jiliang Hu, You He, Matthieu Barbier, Jinyeop Song, and Guy Bunin.
\newblock Transition from global stability to multiple attractors in
  microcosms, 09 2025.

\bibitem[Grilli et~al.(2017)Grilli, Barab{\'a}s, Michalska-Smith, and
  Allesina]{grilli2017feasibility}
Jacopo Grilli, Gy{\"o}rgy Barab{\'a}s, Matthew~J. Michalska-Smith, and Stefano
  Allesina.
\newblock Feasibility and coexistence of large ecological communities.
\newblock \emph{Nature Communications}, 8:\penalty0 14389, 2017.
\newblock \doi{10.1038/ncomms14389}.

\bibitem[Guim and Sonia(2024)]{Aguade2024}
Aguadé-Gorgorió Guim and Kéfi Sonia.
\newblock Alternative cliques of coexisting species in complex ecosystems.
\newblock \emph{Journal of Physics: Complexity}, 5\penalty0 (2):\penalty0
  025022, jun 2024.
\newblock \doi{10.1088/2632-072X/ad506a}.
\newblock URL \url{https://doi.org/10.1088/2632-072X/ad506a}.

\bibitem[Harsanyi(1973)]{Harsanyi1973}
John~C. Harsanyi.
\newblock Oddness of the number of equilibrium points: A new proof.
\newblock \emph{International Journal of Game Theory}, 2\penalty0 (1):\penalty0
  235--250, 1973.
\newblock \doi{10.1007/BF01737572}.

\bibitem[Hill(1973)]{Hill}
M.~O. Hill.
\newblock Diversity and evenness: A unifying notation and its consequences.
\newblock \emph{Ecology}, 54\penalty0 (2):\penalty0 427--432, 1973.
\newblock ISSN 00129658, 19399170.
\newblock URL \url{http://www.jstor.org/stable/1934352}.

\bibitem[Hofbauer and Sigmund(1998)]{HofbaeurSigmund1998}
Josef Hofbauer and Karl Sigmund.
\newblock \emph{Evolutionary Games and Population Dynamics}.
\newblock Cambridge University Press, 1998.

\bibitem[Horn and Johnson(1985)]{Horn_Johnson_1985}
Roger~A. Horn and Charles~R. Johnson.
\newblock \emph{Matrix Analysis}.
\newblock Cambridge University Press, 1985.

\bibitem[Lambert(1992)]{Lambert1992}
Diane Lambert.
\newblock Zero-inflated poisson regression, with an application to defects in
  manufacturing.
\newblock \emph{Technometrics}, 34\penalty0 (1):\penalty0 1--14, 1992.
\newblock \doi{10.2307/1269547}.

\bibitem[Leinster(2024)]{leinster2024entropydiversityaxiomaticapproach}
Tom Leinster.
\newblock Entropy and diversity: The axiomatic approach, 2024.
\newblock URL \url{https://arxiv.org/abs/2012.02113}.

\bibitem[Leinster and Cobbold(2012)]{leincobbo}
Tom Leinster and Christina~A. Cobbold.
\newblock Measuring diversity: the importance of species similarity.
\newblock \emph{Ecology}, 93\penalty0 (3):\penalty0 477--489, 2012.
\newblock \doi{https://doi.org/10.1890/10-2402.1}.
\newblock URL
  \url{https://esajournals.onlinelibrary.wiley.com/doi/abs/10.1890/10-2402.1}.

\bibitem[Lischke and Löffler(2017)]{LISCHKE201724}
Heike Lischke and Thomas~J. Löffler.
\newblock Finding all multiple stable fixpoints of n-species lotka–volterra
  competition models.
\newblock \emph{Theoretical Population Biology}, 115:\penalty0 24--34, 2017.
\newblock ISSN 0040-5809.
\newblock \doi{https://doi.org/10.1016/j.tpb.2017.02.001}.

\bibitem[Lotka(1925)]{Lotka1925}
Alfred~J. Lotka.
\newblock \emph{Elements of Physical Biology}.
\newblock Williams and Wilkins, 1925.

\bibitem[Madec and Gjini(2021)]{MadecGjini2021}
Sten Madec and Erida Gjini.
\newblock The ratio of single to co-colonization is key to complexity in
  interacting systems with multiple strains.
\newblock \emph{Ecology and Evolution}, 11:\penalty0 8456--8474, 2021.
\newblock \doi{10.1002/ece3.7259}.

\bibitem[Magurran(2004)]{Magurran2004}
Anne~E. Magurran.
\newblock \emph{Measuring Biological Diversity}.
\newblock Blackwell Publishing, 2004.

\bibitem[May(1972)]{May1972}
Robert~M. May.
\newblock Will a large complex system be stable?
\newblock \emph{Nature}, 238:\penalty0 413--414, 1972.
\newblock \doi{10.1038/238413a0}.

\bibitem[Mullahy(1986)]{MULLAHY1986341}
John Mullahy.
\newblock Specification and testing of some modified count data models.
\newblock \emph{Journal of Econometrics}, 33\penalty0 (3):\penalty0 341--365,
  1986.
\newblock ISSN 0304-4076.
\newblock \doi{https://doi.org/10.1016/0304-4076(86)90002-3}.

\bibitem[Murty(1972)]{MURTY197265}
Katta~G. Murty.
\newblock On the number of solutions to the complementarity problem and
  spanning properties of complementary cones.
\newblock \emph{Linear Algebra and its Applications}, 5\penalty0 (1):\penalty0
  65--108, 1972.
\newblock ISSN 0024-3795.
\newblock \doi{https://doi.org/10.1016/0024-3795(72)90019-5}.
\newblock URL
  \url{https://www.sciencedirect.com/science/article/pii/0024379572900195}.

\bibitem[Nash(1950)]{Nash1950}
John~F. Nash.
\newblock Equilibrium points in n-person games.
\newblock \emph{Proceedings of the National Academy of Sciences}, 36\penalty0
  (1):\penalty0 48--49, 1950.

\bibitem[Nowak(2006)]{Nowak2006}
Martin~A. Nowak.
\newblock \emph{Evolutionary Dynamics: Exploring the Equations of Life}.
\newblock Harvard University Press, 2006.
\newblock \doi{10.2307/j.ctvjghw98}.

\bibitem[Patil and Rao(1978)]{patil}
G.~P. Patil and C.~R. Rao.
\newblock Weighted distributions and size-biased sampling with applications to
  wildlife populations and human families.
\newblock \emph{Biometrics}, 34\penalty0 (2):\penalty0 179--189, 1978.
\newblock ISSN 0006341X, 15410420.
\newblock URL \url{http://www.jstor.org/stable/2530008}.

\bibitem[Scheffer et~al.(2001)Scheffer, Carpenter, Foley, Folke, and
  Walker]{Scheffer2001}
Marten Scheffer, Steve Carpenter, Jonathan~A. Foley, Carl Folke, and Brian
  Walker.
\newblock Catastrophic shifts in ecosystems.
\newblock \emph{Nature}, 413:\penalty0 591--596, 2001.

\bibitem[Stone(2018)]{Stone1988}
Lewi Stone.
\newblock The feasibility and stability of large complex biological networks: a
  random matrix approach.
\newblock \emph{Sci Rep}, 8, 2018.
\newblock \doi{10.1038/s41598-018-26486-2}.

\bibitem[Strogatz(2018)]{Strogatz2018}
Steven~H. Strogatz.
\newblock \emph{Nonlinear Dynamics and Chaos: With Applications to Physics,
  Biology, Chemistry, and Engineering}.
\newblock CRC Press, 2 edition, 2018.
\newblock ISBN 9780367028760.
\newblock \doi{10.1201/9780429492563}.

\bibitem[Takeuchi(1996)]{Takeuchi}
Y~Takeuchi.
\newblock \emph{Global Dynamical Properties of Lotka-Volterra Systems}.
\newblock WORLD SCIENTIFIC, 1996.
\newblock \doi{10.1142/2942}.
\newblock URL \url{https://www.worldscientific.com/doi/abs/10.1142/2942}.

\bibitem[Taylor and Jonker(1978)]{TaylorJonker1978}
Peter~D. Taylor and Leo~B. Jonker.
\newblock Evolutionary stable strategies and game dynamics.
\newblock \emph{Mathematical Biosciences}, 40\penalty0 (1-2):\penalty0
  145--156, 1978.
\newblock \doi{10.1016/0025-5564(78)90077-9}.

\bibitem[Tilman(1999)]{Tilman1999}
David Tilman.
\newblock The ecological consequences of changes in biodiversity: A search for
  general principles.
\newblock \emph{Ecology}, 80\penalty0 (5):\penalty0 1455--1474, 1999.

\bibitem[Volterra(1926)]{Volterra1926}
Vito Volterra.
\newblock Fluctuations in the abundance of a species considered mathematically.
\newblock \emph{Nature}, 118:\penalty0 558--560, 1926.

\bibitem[Wenger and Freeman(2008)]{WengerFreeman2008}
Seth~J. Wenger and Mary~C. Freeman.
\newblock Estimating species occurrence, abundance, and detection probability
  using zero-inflated distributions.
\newblock \emph{Ecology}, 89\penalty0 (10):\penalty0 2953--2959, 2008.
\newblock \doi{10.1890/07-1127.1}.

\bibitem[Zuur et~al.(2009)Zuur, Ieno, Walker, Saveliev, and Smith]{Zuur2009}
Alain~F. Zuur, Elena~N. Ieno, Neil~J. Walker, Anatoly~A. Saveliev, and
  Graham~M. Smith.
\newblock \emph{Mixed Effects Models and Extensions in Ecology with R}.
\newblock Statistics for Biology and Health. Springer, 2009.
\newblock ISBN 9780387874586.
\newblock \doi{10.1007/978-0-387-87458-6}.

\end{thebibliography}

\end{document}